\documentclass[a4paper]{article}

\usepackage{color,soul}

\usepackage{float}
\usepackage{eurosym}
\usepackage{tcolorbox}
\usepackage{amsmath, amsfonts, amssymb, mathtools}
\usepackage{graphicx,wrapfig}

\usepackage{subcaption}
\usepackage{array}
\newcolumntype{K}[1]{>{\centering\arraybackslash}p{#1}}
\usepackage{dcolumn}
\usepackage{stmaryrd}

\usepackage{adjustbox}
\usepackage{longtable}
\usepackage{empheq}
\usepackage{mathtools}
\usepackage{resizegather}
\usepackage[nogroupskip, nonumberlist, acronym]{glossaries}
\usepackage{amssymb}
\usepackage{tcolorbox}

\title{A machine learning approach to itinerary-level booking prediction in competitive airline markets}
\author{Daniel Hopman, Ger Koole and Rob van der Mei}

\begin{document}

\maketitle

\begin{abstract}
Demand forecasting is extremely important in revenue management. After all, it is one of the inputs to an optimisation method which aim is to maximize revenue. Most, if not all, forecasting methods use historical data to forecast the future, disregarding the "why". In this paper, we combine data from multiple sources, including competitor data, pricing, social media, safety and airline reviews. Next, we study five competitor pricing movements that, we hypothesize, affect customer behavior when presented a set of itineraries. Using real airline data for ten different OD-pairs and by means of Extreme Gradient Boosting, we show that customer behavior can be categorized into price-sensitive, schedule-sensitive and comfort ODs. Through a simulation study, we show that this model produces forecasts that result in higher revenue than traditional, time series forecasts.
\\
\textbf{Keywords:} demand forecasting, effects of competition, traditional statistics, machine learning.
\end{abstract}

\section{Introduction and motivation}

Our definition of RM is the \textit{process of  dynamically assigning capacity to products of perishable nature with a fixed total capacity}. In practice, this means determining what booking classes should be open, for what origin and destination (OD) pair such that the overall, network revenue is maximized. For this optimisation process, we need the airline's own demand forecast, fares and capacity. 

Traditionally, it was thought it was sufficient to segment a market by whether a customer is a business or leisure passenger. An example of this can be found in fare rules: since the inception of pricing, the "Saturday rule" has been used. This rule says that a customer has to stay at least a Saturday night before returning to their origin. Business travellers want to spend their weekends at home, while leisure passengers do not mind spending a Saturday night.

Teichert et al. \cite{teichert2008customer} however, show that this type of segmentation is not sufficient. They find customers that travel in business class for non-business reasons, and customers in economy class travelling for business reasons. Instead, they define five different segments: efficiency, comfort, price, price/performance and all-round performance. Moreover, they find that it is difficult to segment customers, but rather \textit{trips} should be segmented. For example, someone that travels for business may not be price sensitive, but when this same individual travels for leisure, they are. As the industry has changed, whether it is through deregulation, or through the advances in data capture and analytics, segmentation on the other hand, has not. 

Modelling customer behavior is complicated, and there are many reasons for this. First, everyone is different: everyone prioritizes aspects differently. Second, not everyone acts rational. While it is impossible to model these characteristics, it is important to gain an understanding of underlying processes. This can help us understand why people make certain decisions. 

In this paper, we investigate the effects of competition on booking behavior in order to make itinerary-based booking predictions. We do so by combining several data sources. An overview of these data sources is given in the next section. Using these data sources, we engage in feature engineering. Next, we divide these features into those that are airline specific (for example, safety record), while others are itinerary specific (for example, departure time). The objective of this paper is to build a model that given a set of itineraries, to predict what itineraries will be purchased. Armed with this data, airlines can then use this as a strategic tool to increase their demand. 

We consider reviews of a well-known airline review website. This dataset consists of the actual review text, as well as ratings given by the user to the seat, in-flight entertainment (IFE), meal, crew and ground service. Next, we analyze the last $10000$ tweets of the airlines that appear in this article. Based on these two sources, we perform sentiment analysis. 

Another vital dataset gave us access to (historical) pricing information. For a given OD, this dataset captures the price for every airline for every departure date at every day before departure (DBD). Visually, these curves not only tell us what the price was at what point in time before departure, we can also inspect whether airlines react to each other's price change.

We were given access to a data source that includes information about the airline, such as fleet size, fleet age and total aircraft worth. We also have access to an airline safety index.

Lastly, we have OD-specific characteristics. These are features engineered from the OTA's search results. Features include whether this OD has a day flight, whether there is a direct flight, the time of the first departure of the day, the time of the last departure of the day, the number of frequencies, the minimum connection time and the minimum travel time (= flying time + connection time).

This paper is organized as follows. A literature review is conducted in Section \ref{sec:booking_lr}. We provide an overview of the data used for our work in Section \ref{sec:booking_dataoverview}. Our approach to this problem is discussed in Section \ref{sec:booking_approach}, before modelling is covered in Section \ref{sec:booking_mod}. We review the model's performance in Section \ref{sec:booking_results}. A discussion of our work and directions for further research are given in Section \ref{sec:booking_discussion}. In the Appendix, in Section \ref{sec:booking_appendix}, we show all engineered features, its source, and calculation.


\section{Literature review}
\label{sec:booking_lr}

Beckmann \cite{beckmann1958decision} provides the first framework of RM. In this work, published in 1958, he claims that reservations are a form of rationing, and argues that pricing is a direct result of competition. Beckmann, followed by Littlewood \cite{littlewood1972forecasting}, are widely credited as being the first to describe the process of optimisation. In what follows, we separate between findings from surveys and underlying processes, modelling through statistical methods, and modelling through machine learning techniques. 

\subsection*{Surveys and background}
\label{sec:booking_lreview_survey}
In this section, we provide an overview of surveys and the background of the decision making process for bookings. These surveys are conducted in different countries, and aim to illustrate the rationale behind making a booking. The importance of loyalty programmes is also discussed. Other industries, including hospitality are also discussed. Woodside and MacDonald \cite{woodside1994general} provide a framework for the decision and find that the decision of travelling to the destination was made on a separate day than travel plans to that destination were made. This delay in decision making is discussed in \cite{downsell_paper}. The concept of loyalty in the aviation and hospitality industry is very important. The behavior of "loyal" customers is studied by Dolnicar et al.~\cite{dolnicar2011key}. They show that customer satisfaction was not a driver to airline loyalty. They also note that this conclusion only applies for business travellers, and claim it is difficult to model this for leisure travellers. Dowling and Uncles \cite{dowling1997customer} investigate whether loyalty programmes are a driver for success and find any benefits are quickly overshadowed by competition. Sandada and Matibiri \cite{sandada2016investigation} study the service quality and presence of a frequent flyer programme in the airline industry in South Africa and find that customer satisfaction and a frequent flyer programme result in customer loyalty.

The differences between full service, legacy and low cost airlines is discussed by Koklic et al ~\cite{koklic2017investigation}. They show that staff positively influences satisfaction and satisfaction positively influences the likelihood of a customer returning. Tsikriktsis and Heineke \cite{tsikriktsis2004impact} discusses the impact of customer dissatisfaction in the domestic US airline and find evidence that service delivery should be important to an airline, not just wanting to be the cheapest. Park et al. ~\cite{park2005investigating} go one step further and show that that in-flight service, convenience and accessibility were each found to have a positive effect on airline image. Carlsson and Lofgren \cite{carlsson2006airline} study the domestic market of Sweden over a time period of ten years and discuss the cost of a customer churning and show that not every frequent flyer programme is perceived of the same value. The domestic market of India was discussed by Khan et al ~\cite{khan2007customer} and they show that the gap between expectations and perception is of greater importance than the actual product offered. Atalik and Arslan ~\cite{atalik2009study} study the the Turkish market. It was found that flight safety was rated highest, followed by on-time performance, staff and image of the airline. Only after these variables, the price of the ticket was listed in terms of importance. Cho \cite{cho2012impact} studies the domestic United States (US) market. Through surveys, he shows that not only choice of airline, but also choice of airport matters when considering different options. Baker \cite{david2013service} studies the same domestic US market, and investigates the difference between legacy and low cost carriers through data of fourteen US airlines. Baker shows that perceived service quality of low-cost airlines is higher than those of legacy airlines. Keiningham et al. ~\cite{keiningham2014service} discuss the effect of service failures on the future behavior of customers in the domestic US market. They show that major accidents (that cause injuries or fatalities) did not show any impact on customer satisfaction, while minor incidents (such as baggage delays) greatly negatively influence future customer satisfaction. Customer choice in the Nigerian market is studied by Adeola and Adebiyi ~\cite{adeola2014service} and they find that the fare is the most important factor. Khan and Khan \cite{khan2014customer} discuss the market in Pakistan. They show that customer assurance and empathy factors have the largest weights. Hussain et al. ~\cite{hussain2015service} study the United Arab Emirates market and, in particular, how service quality affects customer booking behavior. They show customer expectations has a direct impact on perceived quality, but not customer value. Law \cite{law2017study} investigates the Thai travel market and find that price is the main driver, followed by schedule. The Indonesian market is explored by Manivasugen \cite{manivasugen2013factors}. In this research, he asked 140 students about their travel preferences. The most important factors when deciding between airlines were identified as price, comfort, safety, schedule and airline image.

The growth of the air transport industry as a whole is discussed by Bieger et al. ~\cite{bieger2007driving} It is concluded that ticket price is most important for both economy and  business class passengers. Poelt \cite{poelt2011practical} mentions how traditionally airlines segment customers between leisure and business customers. Teichert et al. ~\cite{teichert2008customer} show why traditional segmentation, both for marketing and optimisation purposes, is no longer sufficient.

The hospitality industry has many similarities with the aviation industry, especially in terms of RM: a customer is given a number of choices, each with a different quality of service. A customer needs to decide what is the best value proposition for them. Verma and Plaschka \cite{verma2003art} show that personalized, on-demand service and brand image are found to be most important. Brunger \cite{brunger2010impact} investigates the "internet price effect" and shows that this effect represents 3 to 8\% in terms of fare value. This reinforces the need to include distribution channel in the definition of RM, which we discuss in \cite{ourpaper_systems}.

Traditionally, it was thought that customers purchasing tickets in premium (first and business class) cabins have a different thought- and decision-making process than customers in economy class cabins. Boetsch et al ~\cite{boetsch2011customer} shows that brand image is more important than both sleep quality (flat seats) and price paid.  

\subsection*{Modelling through statistical methods}
\label{sec:booking_lreview_stat}
In this section, we discuss the modelling of demand forecasts through statistical models. The majority of these methods are based on multinomial logit models.

Coldren et al. ~\cite{coldren2003modeling} study the classification of itinerary through a multinomial logit (MNL) model  in the US domestic market. Vulcano et al. ~\cite{vulcano2010om} study the feasibility of customer choice modelling and the effects of customer choice RM over traditional RM. Lucchesi et al. ~\cite{lucchesi2015airline} evaluate customer preferences for a domestic route in Brazil and try different (logit) methods. Milioti \cite{milioti2015traveler} does not use a MNL approach, but rather a multivariate probit model. Out of those passengers that are not price sensitive, they find that men and business travellers are least likely to be influenced by the level of airfare. Chen et al. ~\cite{chen2011evaluating} argue that the level of service is hard to quantify: words such as "good" or "bad" to describe service quality are ambigious and use fuzzy logic. Ratliff and Gallego ~\cite{ratliff2013estimating} use customer-choice modelling for a different application. They introduce a decision support framework for evaluating sales and profitability impacts of fare brands by using a customer choice model. Airlines only observe bookings if tickets are available for sale. Once the airline stops selling tickets, every demand request is rejected. Estimating this true demand is the process of unconstraining (or, uncensoring).  Haensel and Koole \cite{haensel2011estimating} use a statistical way to estimate unconstrained demand and use real airline data from an airline in The Netherlands. 

\subsection*{Modelling through machine learning}
\label{sec:booking_lreview_ml}
In this section, we give an overview of the research of modelling demand and making predictions through machine learning. 

Mottini and Acuna-Agost ~\cite{mottini2017deep} expose the drawbacks of using a statistical model, specifically the MNL logit model, which was outlined in Section \ref{sec:booking_lreview_stat}. In a similar way, Lhertier et al. ~\cite{lheritier2019airline} discuss itinerary choice modelling. They categorize two kinds of features: features that describe the individual (which they call characteristics) and features that describe alternatives (which they call attributes) and compare MNL with other methods. They show a random forest outperforming a MNL model. Gunarathne et al. ~\cite{gunarathne2018social} show the effects of social media on the service levels of the industry. They find that the number of followers, offensive language used and mention of a competing airline positively influence an airline's response time to an enquiry. Cancellations are modelled through machine learning by Hopman et al. ~\cite{cxl_modelling_paper} They use a combination of traditional statistics and machine learning, in particular extreme gradient boosting, based on features of competition, schedule, price, safety and service quality.

\section{Data overview}
\label{sec:booking_dataoverview}

In this section, we provide a brief overview of the data sources we used.

\subsection{Online Travel Agent (OTA) dataset}
\label{sec:booking_dataoverview_ota}
This data source contains both search queries and bookings made by customers. Customers on this website are tracked through cookies, as well as through their accounts (if they are logged in). In this paper, we only look into actual bookings made, so we omit details of search queries. A sample of the dataset is shown below.

\begin{table}[htbp]
	\centering
	\resizebox{\textwidth}{!}{
	\begin{tabular}{l|l|l|l|l|l|l}
	\hline
		od    & \multicolumn{1}{l}{airline\_id} & \multicolumn{1}{l}{dep\_day\_id} & \multicolumn{1}{l}{dbd} & \multicolumn{1}{l}{dep\_time\_mam} & \multicolumn{1}{l}{travel\_time} & \multicolumn{1}{l}{price} \\
		\hline
		AMS-LHR & 4     & 3063  & -119  & 1305  & 4.33  & 173.92 \\
		AMS-LHR & 4     & 3213  & -3    & 465   & 2.86  & 225.46 \\
		AMS-LHR & 4     & 3444  & -101  & 870   & 7.96  & 178.81 \\
		AMS-LHR & 1     & 3448  & -83   & 420   & 6.82  & 228.74 \\
		AMS-LHR & 2     & 3481  & -33   & 805   & 0.9   & 363.20 \\
		AMS-LHR & 3     & 3621  & -40   & 1265  & 0.9   & 425.37 \\
		AMS-LHR & 4     & 3625  & -98   & 420   & 5.81  & 132.23 \\
		AMS-LHR & 3     & 3677  & -100  & 835   & 0.9   & 453.85 \\
		AMS-LHR & 3     & 3966  & -47   & 865   & 0.9   & 140.14 \\
		AMS-LHR & 2     & 3966  & -47   & 440   & 0.9   & 277.91 \\
		\hline
	\end{tabular}}%
		\caption{Sample of booking dataset}
	\label{tab:sample_booking_data}%
\end{table}%

As we see from Table \ref{tab:sample_booking_data}, we are given itinerary details of bookings. For storage purposes, this company does not store competitor offerings that were not booked.

While the true dataset contains the actual airline name, departure date and other revealing details, this OTA asked to obfuscate airline names and dates for this study, as they consider it sensitive information.

\subsection{Competitive pricing}
The data supplied here contains the historical price for every itinerary offered by every airline. This data is similar in terms of dimensions of the data given in Section \ref{sec:booking_dataoverview_ota}, but contain pricing information of \textit{all} itineraries, not just the one that were purchased. An example of this dataset is given below, in Table \ref{tab:sample_comp_pricing}.

\begin{table}[H]
	\centering
\resizebox{\textwidth}{!}{
	\begin{tabular}{c|c|c|c|c|c|c}
	\hline
		od    & airline\_id & dep\_day\_id & dbd   & dep\_time\_mam & travel\_time & price \\
		\hline
		FRA-SYD & 1     & 946   & -6    & 1220  & 13.24 & 605.73 \\
		FRA-SYD & 2     & 946   & -6    & 1200  & 15.83 & 416.74 \\
		FRA-SYD & 3     & 946   & -6    & 445   & 12.95 & 336.32 \\
		FRA-SYD & 4     & 946   & -6    & 455   & 12.65 & 719.43 \\
		FRA-SYD & 5     & 946   & -6    & 800   & 13.72 & 634.05 \\
		FRA-SYD & 6     & 946   & -6    & 815   & 10.5  & 795.12 \\
		FRA-SYD & 7     & 946   & -6    & 445   & 15.41 & 564.63 \\
		FRA-SYD & 8     & 946   & -6    & 1290  & 14.99 & 677.94 \\
		FRA-SYD & 9     & 946   & -6    & 800   & 14.75 & 582.23 \\
		\hline
	\end{tabular}}%
		\caption{Example of competitive pricing data}
	\label{tab:sample_comp_pricing}%
\end{table}%

\begin{figure}[H]	\includegraphics[width=\textwidth]{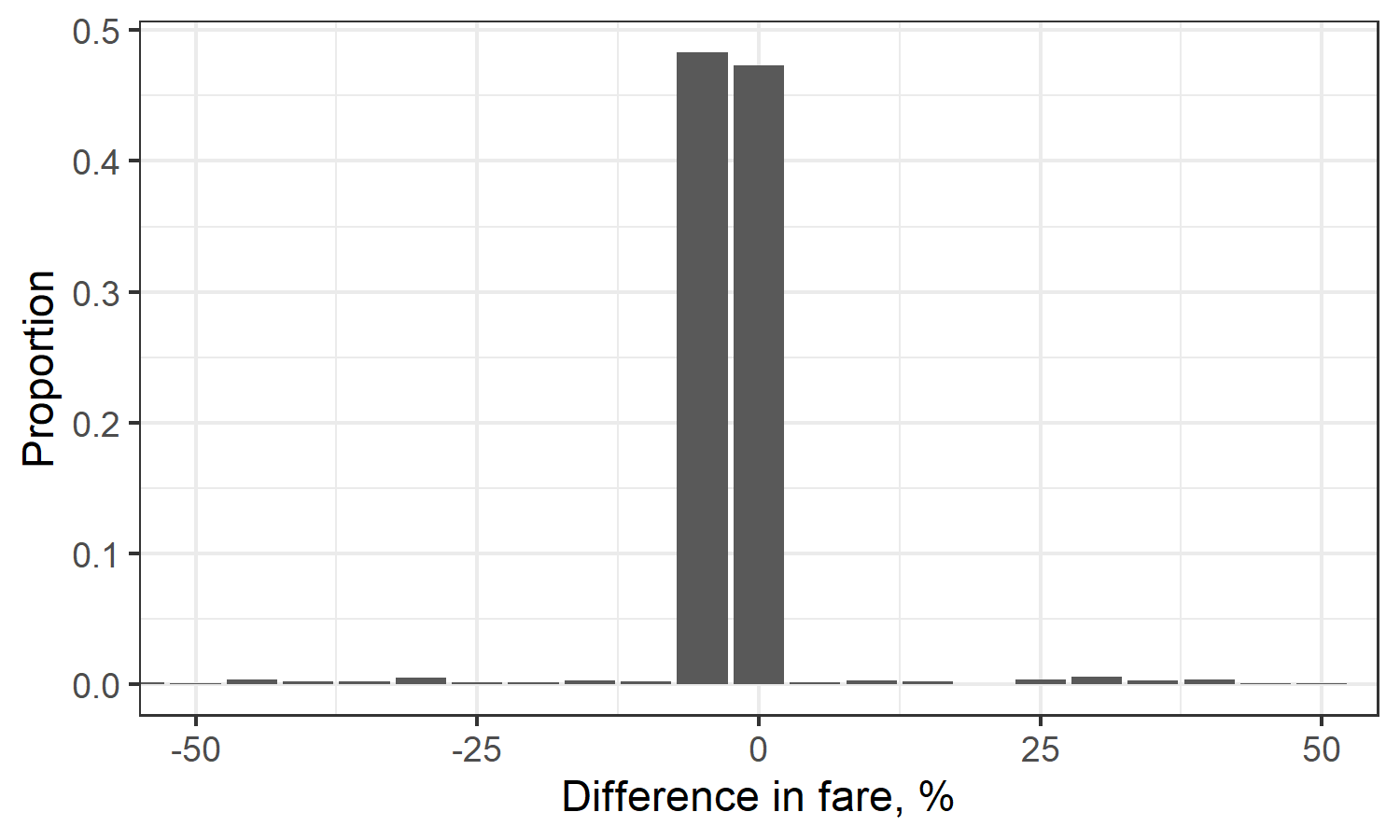} 
	\caption{Fare errors in our dataset}
	\label{fig:fare_errors}
\end{figure}

Note that the data in Table \ref{tab:sample_comp_pricing} greatly enriches the data from the previous section - for every booking made, we now know how each competitor's schedule and price compared. Naturally, the fare in this dataset should, in theory, match with the fare that is associated with the booking. For the vast majority of the bookings ($93\%$ of our dataset is within $1\%$) this is the case : we refer to the reader to Figure \ref{fig:fare_errors}. We review this in our discussion section.

We used the same encoding for obfuscating the airline names and dates. To make comparing fares fair, we always use the fares from this dataset (even though it is given for bookings in our other dataset).

\subsection{Airline ratings}
The website we used for our data contains reviews on airlines. People are given the opportunity to write a (free text) review, as well as rate their trip based on several characteristics, shown below.

\begin{table}[H]
	\centering
\resizebox{\textwidth}{!}{
	\begin{tabular}{c|c|c|c|c|c|c|c|c|c|c}
	\hline
		id & airline & rec. & review & fb    & ground & ife   & crew  & seat  & value & wifi \\
		\hline
	    1     & 5     & N     & (..)  & 2     & 3     & 3     & 2     & 1     & 4     & 1 \\
	    2     & 1     & N     & (..)  & 2     & 2     & 2     & 1     & 1     & 3     & 2 \\
	    3     & 2     & Y     & (..)  & 3     & 4     & 4     & 3     & 4     & 4     & 5 \\
	    4     & 4     & Y     & (..)  & 5     & 4     & 5     & 5     & 5     & 5     & 4 \\
	    5     & 5     & Y     & (..)  & 2     & 2     & 3     & 5     & 4     & 4 & 4 \\
	    \hline
	\end{tabular}}%
		\caption{Airline review dataset. We have omitted the review text to save space. Rec is short for recommended.}
	\label{tab:sample_ratings}%
\end{table}%

Passengers rate their airline based on a general recommendation, their F\&B offering, the service on the ground, the in-flight entertainment (IFE), the quality of the seat, value for money and WiFi performance. They also have the possibility to write a free-text review.

\subsection{Twitter sentiment}
We were given access to the past year's worth of tweets of the airlines present in our data. This dataset required a lot of preparing: we left out retweets (RTs), replies, and we only focused on reviews written in English. Retweets were left out to avoid duplicates records and create a bias towards sentiment. Replies, most often by the airline, are not representative of an individual's perception and for this reason these were left out. Finally, we only focuses on reviews written in English since these would not require translation: translating tweets from other languages may lose the impact they had in their native language.

\subsection{Airline safety ratings}
This website provides an index for airline ratings. It uses accident and incident history, environmental factors and operational risk factors to derive a safety score. A sample of this data is given in Table \ref{tab:sample_safety}.

\begin{table}[htbp]
	\centering

	\begin{tabular}{c|c|c}
	\hline
		rank & airline\_code & score \\
	\hline
		1     & CX    & 0.005 \\
		2     & NZ    & 0.007 \\
		3     & HU    & 0.009 \\
		4     & QR    & 0.009 \\
		5     & KL    & 0.011 \\
		\hline
	\end{tabular}%
		\caption{Airline safety index}
	\label{tab:sample_safety}%
\end{table}%

\subsection{Fleet information}
This dataset contains information on the airline's fleet. It contains several properties of an airline's fleet, such as size and cost. While not always correct, we aim to use the average fleet's age as a proxy for a comfort rating (newer aircraft are typically quieter and provide better entertainment). Similarly, we intend to use the fleet size as a proxy to how well passengers are accommodated when irregular operations happen (if an airline only has a handful of aircraft and a flight gets cancelled, it is likely a passenger will endure long delays). An example of this data is shown in Table \ref{tab:sample_airlineinfo}.

\begin{table}[H]
	\centering
\resizebox{\textwidth}{!}{
	\begin{tabular}{c|c|c|c|c}
	\hline
		\multicolumn{1}{c}{airline\_id} & aircraft & \multicolumn{1}{c}{aircraft cost} & aircraft registration & \multicolumn{1}{c}{aircraft age} \\
		\hline
		1     & 77W   & 300   & PH-ABC & 8.8 \\
		1     & 77W   & 300   & PH-XYZ & 12.2 \\
		\hline
	\end{tabular}}
		\caption{Fleet information example}
	\label{tab:sample_airlineinfo}%
\end{table}%

In Table \ref{tab:sample_airlineinfo}, we show the airline ID, again, obfuscated to anonymize the data, the IATA aircraft code (for example, 77W represents a Boeing 777-300ER), its obfuscated registration, and aircraft age in years.

\subsection{Data overview}

In this section, we provide a few characteristics of our dataset. Table \ref{tab:num_comp} shows the number of competitors by OD. This is not an exhaustive list of all airlines that sell this OD. Rather, it is a list of airlines that sold this OD in our dataset. 

\begin{table}[H]
\centering
\resizebox{\textwidth}{!}{
\begin{tabular}{cccccccccc}
  \hline
AMS-DXB & AMS-LHR & AMS-SYD & CDG-SYD & FRA-SYD & FRA-KUL & FRA-SYD & KUL-SIN & LHR-JFK & LHR-SYD \\ 
 7 & 4 & 5 & 4 & 9 & 6 & 5 & 2 & 2 & 5 \\ 
   \hline
\end{tabular}}
\caption{Number of competitors by OD}
	\label{tab:num_comp}%
\end{table}

\begin{figure}[H]
	\includegraphics[width=\textwidth]{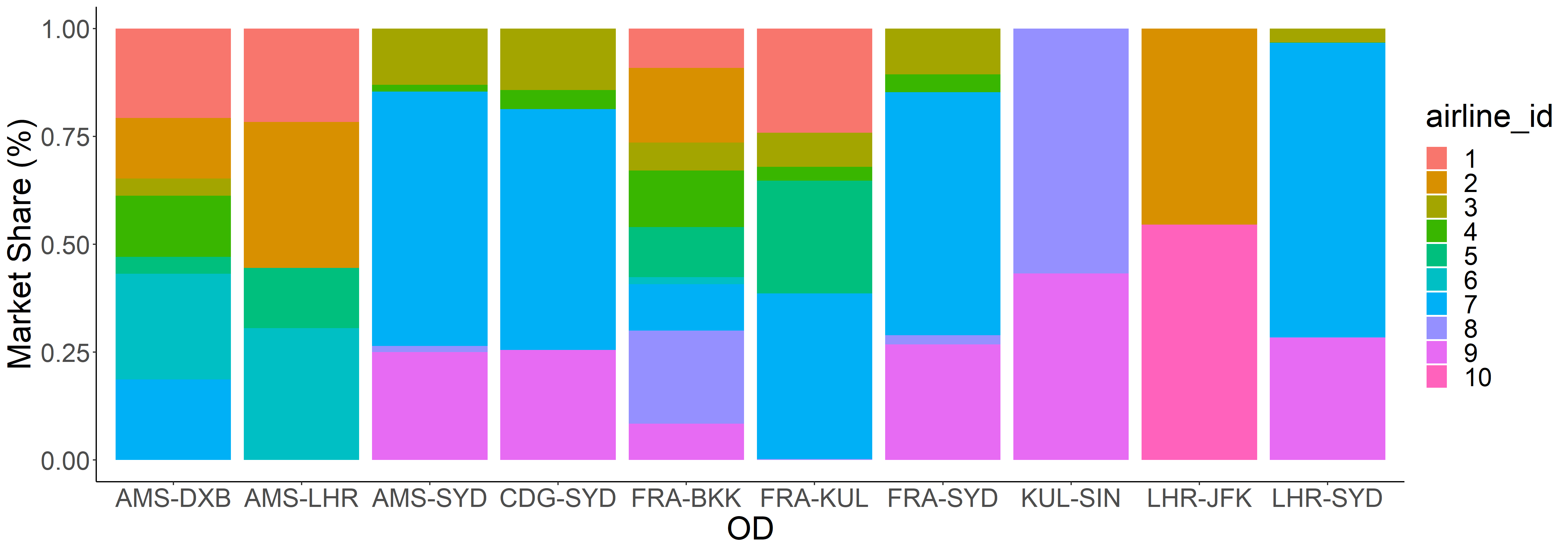} 
	\caption{Demand by OD and competitor}
	\label{fig:demand_by_comp}
\end{figure}

Figure \ref{fig:demand_by_comp}  shows the demand (sum of bookings) of competitors by OD. Note not every airline is operating every OD. For example, it is unlikely that an airline operating KUL-SIN will also operate LHR-JFK. 

\begin{figure}[H]
	\includegraphics[width=\textwidth]{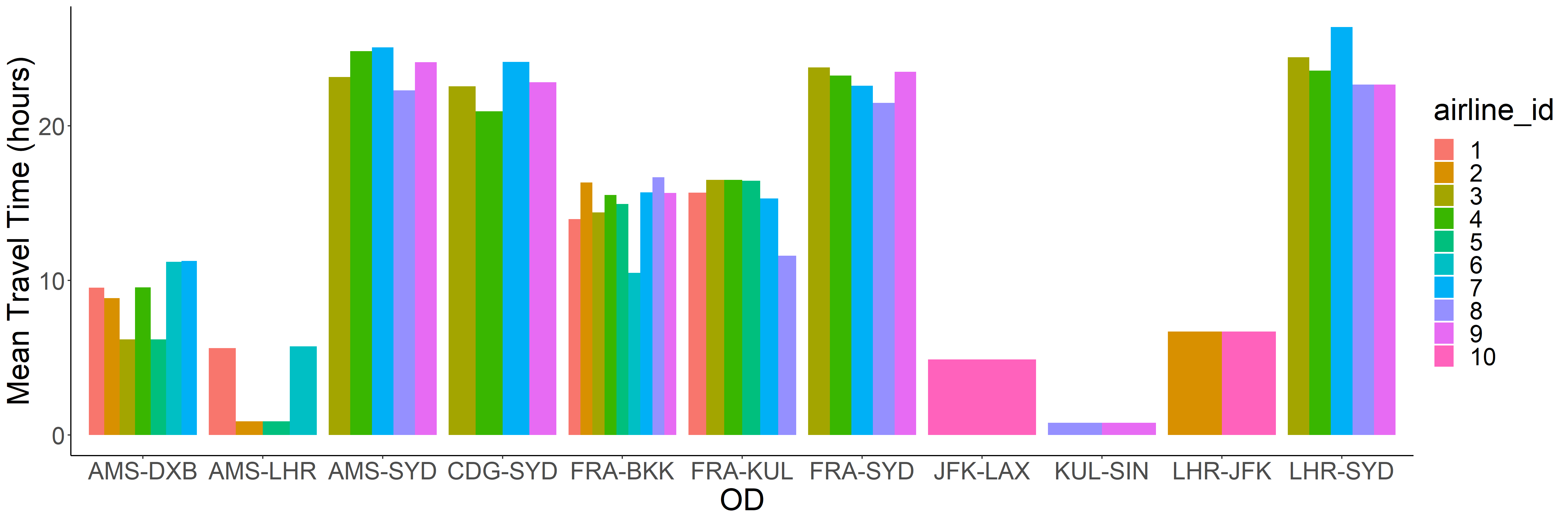} 
	\caption{Travel time by OD and competitor}
	\label{fig:traveltime}
\end{figure}

Figure \ref{fig:traveltime} illustrates the average travel time in hours by airline by OD. 

Taking a look at AMS-LHR, the big difference can be explained by the fact that airline 2 and 3 offer a direct flight between these airports, while airlines 1 and 4 offer connections, resulting in longer average travel times. Note that for LHR-SYD the small difference in travel times across airlines. Even though there are five competing airlines in this OD, each with their own hub, and a thus a different way to fly, the resulting travel time is quite similar.

\section{Approach}
\label{sec:booking_approach}
In this section, we will illustrate the approach we have taken. First, we reiterate the objective in Section \ref{sec:booking_approach_objective}. Next, in Section \ref{sec:booking_approach_featureng} we discuss the features we have engineered.

\subsection{Objective}
\label{sec:booking_approach_objective}
Consider an itinerary, with features as described in the previous section. The objective is to determine whether an itinerary will be purchased. We will use Extreme Gradient Boosting to build this binary logistic model (after all, an itinerary is purchased yes or no). Note that this is different from typical choice modelling. In that case, the objective is given a set of alternatives, what option a customer is most likely to purchase. Our objective goes beyond that - we look at the features of an individual itinerary to determine whether this will be purchased.

Since we suspect that behavior will be different across ODs, we suspect we cannot compare a short-haul to a long-haul OD, or a business-heavy to a leisure-heavy OD, we build a separate model for every OD.

\subsection{Feature engineering}
\label{sec:booking_approach_featureng}
Arguably, the most important step in building the model is the process of feature engineering. Feature engineering "involves constructing novel features from given data (...) driven by domain knowledge developed over time (...) with the goal of improving predictive learning performance." ~\cite{khurana2016cognito}. In this section, we follow the same structure as in the Data Overview section, Section \ref{sec:booking_dataoverview}, and will review our feature engineering.

\subsubsection{OTA dataset}

From this dataset, we construct the following \textbf{airline/OD-specific features}:

\begin{enumerate}
	\item Number of unique itineraries offered
	\item Minimum, maximum and average flying time
	\item Minimum, maximum and average travel time (flying time+connection time)
	\item Is an evening departure offered (boolean, if there is an itinerary departure time after 6PM)
	\item Departure time of first flight of the day (first sector), in minutes after midnight
	\item Departure time of last flight of the day (first sector), in minutes after midnight
	\item Arrival time of first flight of the day (first sector), in minutes after midnight
	\item Arrival time of last flight of the day (first sector), in minutes after midnight
	\item Is a direct flight offered? (boolean, true if and only if there is an itinerary with one sector)
\end{enumerate}

Next, we derive the following \textbf{itinerary-specific features}:

\begin{enumerate}
	\item Flying time, in hours
	\item Travel time, in hours
	\item Morning, afternoon or evening departure time: departure time for the first sector before 9AM, before 6PM, or after 6PM respectively
	\item Morning, afternoon or evening arrival time: arrival time for the first sector before 9AM, before 6PM, or after 6PM respectively
	\item Wide-body or narrow-body aircraft used
	\item Night flight (boolean, if departure time is before midnight and arrival time is after midnight, both in local timezones).
\end{enumerate}

\subsubsection{Competitive pricing}

Using this dataset, we derived features that measure current and rolling window price movements. In the example that follows, we look into the calculation of pricing features for airline 1 for 100 days before departure. The features compare this airline's pricing against the cheapest airline in the market, dubbed $yy$, and second cheapest airline in the market, dubbed $xx$. We use this terminology, since airline codes typically consist of two characters.

Table \ref{tab:pricing_features_1fares} shows a subset of fares for four different airlines at specific times before departure. 

\begin{table}[H]
	\centering

	\begin{tabular}{c|c|c|c|c}
		\hline
		time  & airline\_1 fare & airline\_2 fare & airline\_3 fare & airline\_4 fare \\
		\hline
		-103  & 450   & 600   & 500   & 1000 \\
		-102  & 475   & 600   & 500   & 1100 \\
		-101  & 450   & 625   & 360   & 1100 \\
		-100  & 450   & 750   & 500   & 300 \\
		\hline
	\end{tabular}%
	\caption{Example of fares}
	\label{tab:pricing_features_1fares}%
\end{table}%

\begin{table}[H]
	\centering
	  \resizebox{\textwidth}{!}{
	\begin{tabular}{c|c|c|c|c|c|c}
	\hline
		time  & airline\_1 fare & yy fare & xx fare & is\_cheapest? & yy fare difference & xx fare difference \\
		\hline
		-103  & 450   & 450 (1) & 500 (3) & Y     & 0     & -50 \\
		-102  & 475   & 475 (1) & 500 (3) & Y     & 0     & -25 \\
		-101  & 450   & 360 (3) & 450 (1) & N     & 90    & 0 \\
		-100  & 450   & 300 (4) & 450 (1) & N     & 100   & 0 \\
		\hline
	\end{tabular}}%
	\caption{Step 1 - Calculation of fare difference. The number between parentheses indicates what airline this fare belongs to.}
	\label{tab:pricing_features_2calc}%
\end{table}%

Table \ref{tab:pricing_features_2calc} shows how we first calculate the $yy$ and $xx$ fare. This is simply the cheapest and second cheapest fare. Note that the parentheses indicate what airline offers that specific fare. Our first engineered feature is \textit{is\_cheapest}, which is a Boolean and is true if the airline's fare equals the $yy$ fare. Note that you can have multiple airlines that are cheapest, if they have the same fare. Next, we calculate the difference between the airline and the $yy$ and $xx$ fares.

\begin{table}[H]
	\centering
	
	\begin{tabular}{c|c|c|c}
	\hline
		time  & mean3d\_yy & mean\_3d\_xx & farediff \\
		\hline
		-100  & (0+0+90)/3=30 & (-50-25+0)/3=-25 & 450-360=90 \\
		\hline
	\end{tabular}%
	\caption{Step 2 - Feature engineering}
	\label{tab:pricing_features_3featureeng}%
\end{table}%

Table \ref{tab:pricing_features_3featureeng} then illustrates how these features are calculated. Note that the $mean3d_{yy}$ and $mean3d_{xx}$ fares indicate the average difference in fares over the past three days, while $farediff$ is a snapshot feature that only measures the difference at only one particular day ($t=-100$).

Apart from calculating the mean, we also calculate the standard deviation. We repeat this process for a $7$, $14$ and $28$ rolling window. As a result, we will have the following features:

\begin{enumerate}
	\item $mean3d_{yy}, sd3d_{yy}, mean3d_{xx}, sd3d_{xx}$
\item $mean7d_{yy}, sd7d_{yy}, mean7d_{xx}, sd7d_{xx}$
\item $mean14d_{yy}, sd14d_{yy}, mean14d_{xx}, sd14d_{xx}$
\item $mean28d_{yy}, sd28d_{yy}, mean28d_{xx}, sd28d_{xx}$
\item \textit{farediff}.
\item \textit{is\_cheapest}
	
\end{enumerate}

\subsubsection{Airline ratings}

From our dataset, we derived median values for the characteristics, as well as counted the observations. These are given in Table \ref{tab:airline_ratings}.

\begin{table}[htbp]
	\centering
	\resizebox{\textwidth}{!}{
	\begin{tabular}{c|c|c|c|c|c|c|c|c|c|c}
	\hline
		airline & rec. & review & f\&b    & ground & ife   & crew  & seat  & value & wifi  & obs \\
		\hline
		1     & 0.89  & 6.52  & 4.19  & 4.18  & 3.6   & 4.65  & 3.8   & 4.1   & 2.6   & 79 \\
		2     & 0.52  & 5.84  & 3.22  & 3.04  & 2.79  & 3.64  & 3.25  & 3.25  & 2.66  & 638 \\
		3     & 0.74  & 6.21  & 3.62  & 3.77  & 4.44  & 3.73  & 3.97  & 3.75  & 3.63  & 221 \\
		4     & 0.6   & 5.63  & 3.15  & 2.94  & 3.29  & 3.37  & 3.31  & 3.06  & 3.2   & 87 \\
		5     & 0.89  & 6.35  & 3.84  & 3.89  & 3.35  & 4.19  & 3.86  & 4.2   & 2.67  & 71 \\
		\hline
	\end{tabular}}%
	\caption{Example of airline ratings}
	\label{tab:airline_ratings}%
\end{table}%

Some reviews contain what route the reviewer flew. Ideally, one would only look into reviews that match the OD we are studying. However, the resulting number of reviews are too low to be reliable measures for aggregates. For this reason, we have chosen to take aggregates across airlines. Note in Table \ref{tab:airline_ratings} that there is no OD present. We review this decision in our discussion.

Table \ref{tab:airline_ratings} shows whether passengers recommend this airline, as well as the median scores for the onboard F\&B, the service provided on the ground, the IFE, the crew, quality of the seat, value for money and WiFi. The "obs" column show how many reviews we collected. 

The value under review is constructed using sentiment analysis. We do this as follows. The free text of all reviews are read into $R$. Each review is converted into a long $1 * N$ vector by splitting the review; each element in the vector will have a single word. First, stop words and punctuation are removed. The text is converted to lowercase. This vector may then be joined with the AFINN dataset. The AFINN dataset \cite{nielsen2011new} was created by Nielsen containing 1468 unique words, that Nielsen manually labeled with a score between minus five (highly negative) and plus five (highly positive). Next, we simply take the mean over this list of scores to determine how positive or negative a review was, and scale this value between $0$ and $10$. An example of the AFINN dataset is shown in Table \ref{tab:sample_afinn}.

\begin{table}[htbp]
	\centering
	
	\begin{tabular}{l|c}
	\hline
		word  & \multicolumn{1}{c}{score} \\
		\hline
		amazing & 4 \\
		breathtaking & 5 \\
		disaster & -2 \\
		distrust & -3 \\
		excellence & 3 \\
		fraudsters & -4 \\
		limited & -1 \\
		misleading & -3 \\
		\hline
	\end{tabular}%
	\caption{AFINN subset sample}
	\label{tab:sample_afinn}%
\end{table}%

\subsubsection{Twitter sentiment}

After removing quotes, retweets, special characters among other things, we followed the same procedure as for airline ratings reviews: we matched every word in the tweets with the AFINN list, then took an average of these scores to get a rating by tweet. We aggregated these individual reviews by taking the median of each tweet's rating to obtain a score by airline, then scaled it to a value between $0$ and $10$. Just like for airline ratings, we were unable to obtain scores by OD - we were unable to derive from the tweets what OD passengers were flying. An example of this dataset is given in Table \ref{tab:sample_twitter}.

\begin{table}[htbp]
	\centering
	
	\begin{tabular}{c|c}
	\hline
		\multicolumn{1}{l}{airline\_id} & \multicolumn{1}{l}{twitter\_sentiment} \\
		\hline
		1     & 6.32 \\
		2     & 4.22 \\
		3     & 7.02 \\
		4     & 6.49 \\
		5     & 7.21 \\
		\hline
	\end{tabular}%
	\caption{Twitter sentiment scores by airline}
	\label{tab:sample_twitter}%
\end{table}%

\subsubsection{Airline safety Ratings}

We did not engineer any features. Instead, we used the score provided to us.

\subsubsection{Fleet information}

For fleet information, we derived the following \textbf{airline-specific features}:

\begin{enumerate}
	\item Fleet size: number of aircraft
	\item Fleet cost: sum of aircraft cost (list price)
	\item Fleet age: median of aircraft age
\end{enumerate}

An example of these features is shown in \ref{tab:featuresfleet}.

\begin{table}[htbp]
	\centering
	\begin{tabular}{c|c|c|c}
	\hline
		airline\_id & fleet size & fleet cost & fleet median age \\
		\hline
		a & 226   & 42343 & 11.19 \\
		b & 268   & 54130 & 11.67 \\
		c & 249   & 85298 & 6.15 \\
		\hline
	\end{tabular}%
	\caption{Fleet information example}
	\label{tab:featuresfleet}%
	
\end{table}%

In Table \ref{tab:featuresfleet} we have used a different $airline_{id}$ as in previous sections, as there are ways to trace the actual airline's name using these characteristics.

\section{Modelling}
\label{sec:booking_mod}
Having engineered features, we will predict whether a given itinerary is purchased. To accomplish this, we will use extreme gradient boosting.

\subsection{Extreme gradient boosting}
\label{sec:chapter_booking_behavior_xgb}
In this section, we provide a brief overview of extreme gradient boosting (XGB). For a full introduction of XGB we refer the reader to Chen et al \cite{chen2016xgboost}. In what follows, we provide a short, alternative brief. Suppose we have an input $x_i$, and an output $y_i$. We would like to make a prediction, denote this by $\widehat{y_i}$. 

An example of estimating $y_i$ is given in Equation (\ref{eq:booking_regression}).
\begin{equation}
    \widehat{y_i} = \sum_{j=1}^{} \alpha_j x_{ij}
    \label{eq:booking_regression}
\end{equation}
This, of course, is simple linear regression (we omit a constant and standard error): in this case $y_i$ is expressed as a linear combination of explanatory variables, denoted by $x_i$. The objective is to estimate those $\alpha_j$ that minimize an error measure. We typically want to minimize some error measure, denoted by $O_m$, depending on a model $m$. A natural selection of an objective function is an error measure:

\begin{equation}
    O_m(\alpha) = \sqrt{\sum_{j=1}^{}(y_i(\alpha) - \widehat{y_i})^2}.
    \label{eq:booking_errormeasure}
\end{equation}

This is known as the root mean squared error (RMSE).

Extreme Gradient Boosting is a tree boosting algorithm. The method works in a similar fashion: first, specify how our predictor is expressed in terms of features, like
(\ref{eq:booking_regression}). Next, specify an objective function, comparable to (\ref{eq:booking_errormeasure}). Finally, iterate to find the optimal value.

Consider the following objective function:

\begin{equation}
    O_{XGB}(\theta) = L(\theta) + \Omega(\theta)
    \label{eq:booking_obj_xgb}
\end{equation}

In Equation (\ref{eq:booking_obj_xgb}), the first term, $L(\theta)$ is the loss function. Typically, (R)MSE is used. The second term, $\Omega(\theta)$, is the regularization term. This term measures the complexity of the model and helps us control overfitting.

Let $f_k(x_i)$ be a function that takes a set of input parameters $x_i$, as before, and return the score of $x_i$ in tree $k$, $k = 1, .. , K$. Suppose we have $M$ different features. In XGB, we assume that
\begin{equation}
    \widehat{y_i}  = \sum_{k=1}^{K} f_k(x_i).
    \label{eq:booking_xgb_yi}
\end{equation}

That is, the prediction for $y_i$, $\widehat{y_i}$, is the sum of linear combination of the score in each tree. Let $T$ be the number of leaves in a tree, and $w_i$ the weight of leaf $i$. Assume that we have $M$ features. Therefore, our input $x_i$ is a M-dimensional vector. Introduce a function $q(x_i)$ which takes an input $x_i$, and follow the rules of the decision tree to map it to the leaves. Specifically:
\begin{equation*}
q(x_i), \mathbb{R}^M \longrightarrow T
\end{equation*}
The prediction is then given by this function $q(x_i)$, weighted by the weights of the leaves, denoted by $w_i$. Therefore:
\begin{equation*}
f(x_i) = q(x_i)w_i
\end{equation*}
Since XGB is an iterative method, there is no static objective function as in Equation (\ref{eq:booking_errormeasure}). Similarly, the prediction of $y_i$ at time $t$ is given by the previous value of $y$, represented as $y_{i}^{t}$, plus the score of $x_i$ in our new tree:
\begin{equation}
    \widehat{y_{i}^t}  = \widehat{y_{i}^{t-1}} + f_t(x_i)
    \label{eq:booking_xgb_yi_iter}
\end{equation}
Suppose now we have a generic loss function, $l$, some choice for $L$, as we introduced in (\ref{eq:booking_obj_xgb}). In this case, we have:
\begin{equation}
   O_{xgb}^t  = \sum_{i=1}^{n} \Big( l  ( y_i, \widehat{y_{i}^{t-1}} + f_t(x_i) ) \Big) + \Omega(f_k)
    \label{eq:booking_xgb_obj_generic_t}
\end{equation}

XGB uses a second-order Taylor expansion to approximate this function $l$. Recall that the Taylor expansion of $f(x)$ at $x+a$ up to the second degree is given by:
\begin{equation}
 f(x) = f(x+a) + \frac{f'(x+a)}{1!}(x+a) + \frac{f''(x+a)}{2!}(x+a)^2
    \label{eq:booking_xgb_taylor}
\end{equation}

Deciding to use the MSE for our generic function $l$, our objective function $O_{XGB}$ at time t is equal to:
\begin{align}
   O_{XGB}^t  = &  \sum_{i=1}^{n} \Big( y_i -  ( \widehat{y_{i}^{t-1}} + f_t(x_i) ) \Big)^{2} + \Omega(f_k) \\
               = & \sum_{i=1}^{n} \Big(2 (\widehat{y_{i}^{t-1}} - y_i) f_t(x_i) + f_t(x_i)^2 \Big) + \Omega(f_k).
    \label{eq:booking_xgb_obj_t}
\end{align}

The term that remains is $\Omega(f_k)$. As we discussed above, this term is important but often forgotten, and helps us control the complexity of the models by penalizing large models. In XGB, this function is defined as follows:
\begin{equation}
   \Omega(f_k) = \gamma T + \dfrac{1}{2}\lambda \sum_{i=1}^{n} (w_i)^2.
    \label{eq:booking_xgb_omega}
\end{equation}

In the regularization shown in (\ref{eq:booking_xgb_omega}), $\gamma$ is threshold of reduction in the loss function for XGB to further split a leaf. Smaller values will make XGB split more leaves, therefore generating a more complex tree structure, while larger values will limit the number of leaves. We chose $\gamma = 0.25$. On the other hand, $\lambda$ penalizes on large values of $w_i$. Intuitively, this is an appealing property: it encourages XGB to use all of its inputs a little bit, rather that some of its inputs a lot. The choice of $\lambda$ is defined by the user. In our work, we chose $\lambda = 1$. Other parameters are investigated in Section \ref{sec:booking_select_params}.

\subsection{Selecting parameters}
\label{sec:booking_select_params}

\textbf{Learning Rate, $\eta$}\\
The learning rate is the shrinkage used. The shrinkage factor is a way to slow down the incremental performance gain of a new tree being added to the ensemble. A smaller learning rate means the model will take longer to run but is less likely to overfit. We try values of $\eta = 0.01, 0.02 , .. , 0.1$.
\\

\noindent\textbf{Number of Decision Trees, $n_t$}\\
The number of decision trees specifies how many trees can be used until we stop optimizing. In practice, this number is typically relatively low, $1000$ or less (depending on the model size) and is a direct result of how the algorithm works. More specifically, it is the result of how fast errors are being corrected. A new boosted tree model is constructed based on errors of the current tree. We therefore expect to see diminishing returns. Let $n_t$ be the number of trees we can use. We perform a grid search over values of $n_t = 50, 100, .. , 500$.\\

\noindent\textbf{Depth of the tree, $d_t$}\\
This parameter specifies how many layers a tree may have. Intuitively, a small number of layers in a tree do not capture enough details about the data to be a good descriptor. On the other hand, a tree with too many levels may be overfitting the dataset. Let $d_t$ be the depth of the tree. We will evaluate values of $d_t = 3, 4, ... , 20$.\\

\noindent\textbf{Subsample, $s_t$}\\
Subsample represents the percentage of the number of observations chosen to fit a tree. Taking out too much data means the model will run faster (after all, less fitting needs to be done), but taking not enough data may expose us to overfitting. Let $s_t$ be the proportion of data used to fit a tree. Then we will try $s_t = 0.2, 0.3, ... , 1$.\\

\noindent\textbf{Number of features used per tree, $f_t$}\\
In the dataset, each row represent an observation. Every column contains a feature. The XGB algorithm samples the number of columns in building a new tree.  Using all columns for every tree may lead to overfitting, but also makes the problem slower to solve.\\

We present the results in Section \ref{sec:booking_results}. More specifically, in Section \ref{sec:booking_xgboost_perf}, we will review the influence of each of these parameters for the different ODs. 

\subsection{Performance}
\label{sec:booking_xgb_performance}

Suppose we have the following, generalized linear model that is used to make a prediction for a value $y_i$:

\begin{equation}
\hat{y_i} = \sum_{j}{}{w_j x_{ij}}
\end{equation}

In our dataset, we have an exhaustive list of options offered by the most popular airlines on this route. We know what itinerary was purchased. This will be our label. As a result, we have a logistic binary objective function: yes or no. 

The extreme gradient boosting algorithm will return a probability of purchase – between 0 and 1. To obtain a yes/no label, we round the probability to the nearest integer. We used the R implementation of the widely used \textit{xgboost} package, specifically version $0.72$. We used a maximum tree depth of $20$, and used $10$ passes (iterations) of the data.

\begin{figure}[H]
    \begin{center}
	\includegraphics[width=0.9\textwidth]{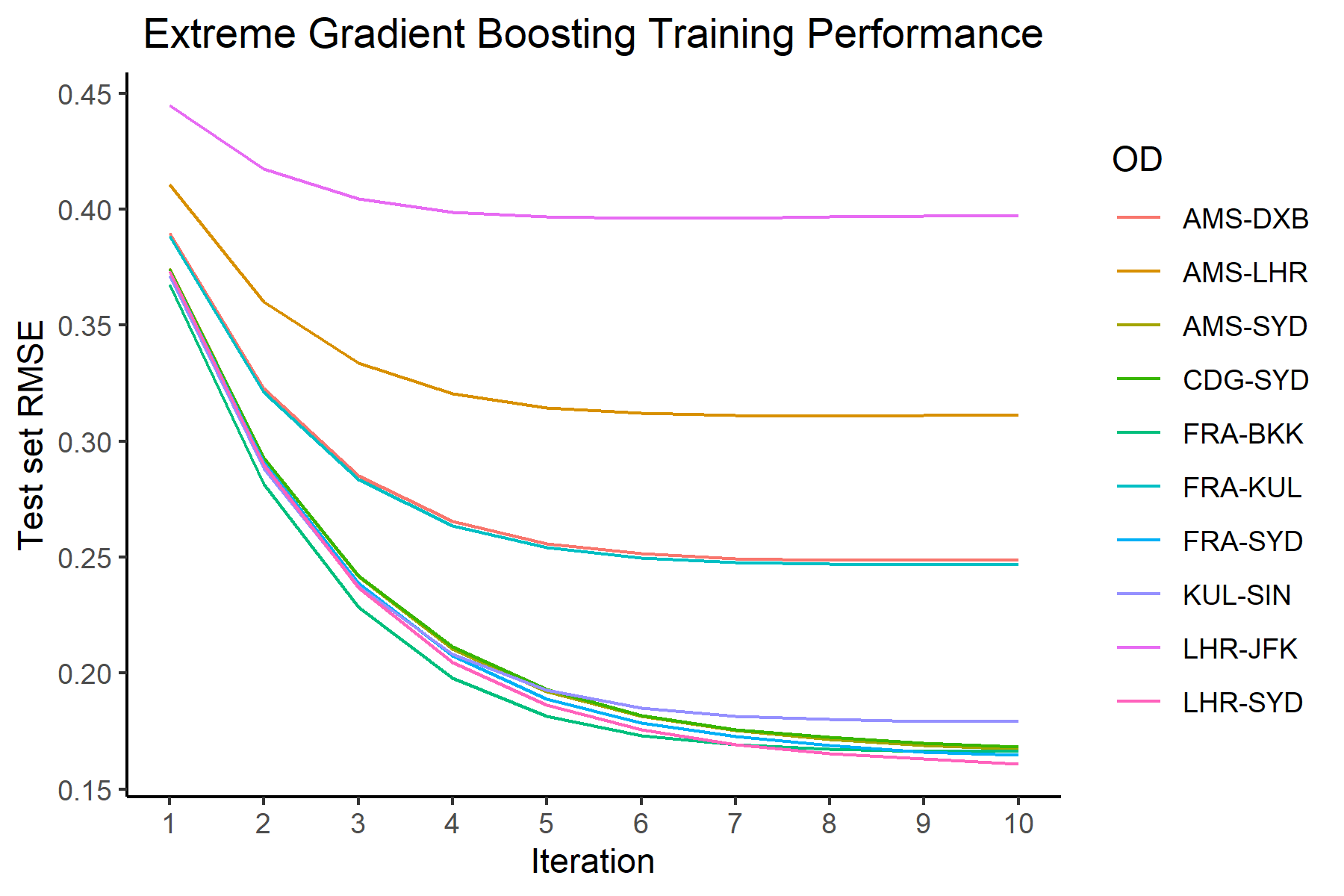} 
	\end{center}
	\caption{RMSE error improvement by iteration}
	\label{fig:error_interation}
\end{figure}

Figure \ref{fig:error_interation} shows the RMSE by iteration. Note from this graph that the RMSE decreases quickly. We show that the number of passes chosen, $10$ is sufficient to obtain sufficiently close to optimal performance, while making sure runtime is acceptable.

\subsection{Runtime}
The XGB model performance is excellent in terms of runtime - we record an average of $11.9$ seconds across the 10 ODs we consider.

On the other hand, the feature engineering process takes a significantly longer amount of time. This involves a series of expensive operations on the dataset. Calculating the pricing features, in particular, involves joining a dataset on itself. This was done in memory in R. Despite this, the entire process, from loading data to the actual engineering of features took an average of $23$ minutes per OD. This is important metric, since frequent forecasting and reoptimisation (discussed in \cite{ourpaper_lit}) is important in practice.

\section{Results}
\label{sec:booking_results}

\subsection{Comparison with logistic regression}
In this section, we compare the performance of the proposed model against the logistic regression model. The results are shown in Table \ref{tab:booking_confmat_logit}. Note that the model is biased against predicting a non-purchase, since the majority of itineraries are never purchased. We therefore omit true negatives from our performance metric, and study the proportion of false negative, false positive, and true positive predictions.

\begin{table}[H]
\centering
\begin{tabular}{l|rr|rr|rr}
  \hline
 Classification & \multicolumn{2}{c}{False Negative} & \multicolumn{2}{c}{False Positive} & \multicolumn{2}{c}{True Positive} \\
  OD / Method  & Logit & XGB & Logit & XGB & Logit & XGB \\
  \hline
 AMS-DXB & 0.55 & \textbf{0.34} & \textbf{0.05} & 0.19 & 0.39 & \textbf{0.47} \\ 
 AMS-LHR & 0.40 & \textbf{0.26} & \textbf{0.10} & 0.18 & 0.50 & \textbf{0.56} \\ 
 AMS-SYD & 0.63 & \textbf{0.11} & \textbf{0.05} & 0.06 & 0.32 & \textbf{0.83} \\ 
 CDG-SYD & 0.39 & \textbf{0.09} & 0.10 & \textbf{0.06} & 0.51 & \textbf{0.86} \\ 
 FRA-SYD & 0.24 & \textbf{0.17} & \textbf{0.07} & 0.13 & 0.69 & \textbf{0.70} \\ 
 FRA-KUL & 0.36 & \textbf{0.23} & \textbf{0.11} & 0.19 & 0.53 & \textbf{0.58} \\ 
 FRA-SYD & 0.80 & \textbf{0.10} & \textbf{0.04} & 0.06 & 0.16 & \textbf{0.83} \\ 
 KUL-SIN & 0.09 & \textbf{0.07} & \textbf{0.00} & 0.02 & 0.90 & \textbf{0.91} \\ 
 LHR-JFK & 0.39 & \textbf{0.20} & \textbf{0.10} & 0.17 & 0.51 & \textbf{0.62} \\ 
 LHR-SYD & 0.63 & \textbf{0.10} & \textbf{0.05} & 0.06 & 0.32 & \textbf{0.84} \\ 
   \hline
\end{tabular}
\caption{Comparison of False Negatives, False Positive and True Positives for the logistic regression (logit) and XGB model. Bolded values represent the better value}
\label{tab:booking_confmat_logit}
\end{table}

Table \ref{tab:booking_confmat_logit} shows the performance by OD. The numbers in bold face compare the logit and XGB model by prediction type and highlight the better value. First, let us consider the false negative: the XGB model outperforms the logit model for every OD. The differences range from $0.02$ for KUL-SIN to $0.7$ for FRA-SYD. We suspect that the logit model for the KUL-SIN OD performs relatively well since the number of competitors is low (namely, 2) and seems to be driven by a single feature, which we will discuss in Section \ref{sec:booking_customerbehavior}.

Reviewing the false positives, there is only one OD in which the number of false positives are lower for the XGB model. The differences range from outperforming logit by 4\% to 14\% more false positive predictions for AMS-DXB. Comparing these false positives to false negatives, the results seem to indicate that the logitistic regression model is biased toward predicting false negatives, while the XGB model is biased toward predicting false positives. 

However, most interesting are the true positives. In all cases, the XGB model outperforms the logit model. This seems to illustrate the need for a more advanced method than simple logistic regression. The differences in performance range from 1\% on the KUL-SIN OD, to 67\% on FRA-SYD. For FRA-SYD, note that the number of false negatives is almost equivalent to the number of true positives. In fact, the performance gains on most SYD ODs are impressive: 51\% for AMS-SYD, 15\% for CDG-SYD and 52\% for LHR-SYD.

\subsection{Overall performance}

\begin{figure}[H]
	\includegraphics[width=\textwidth]{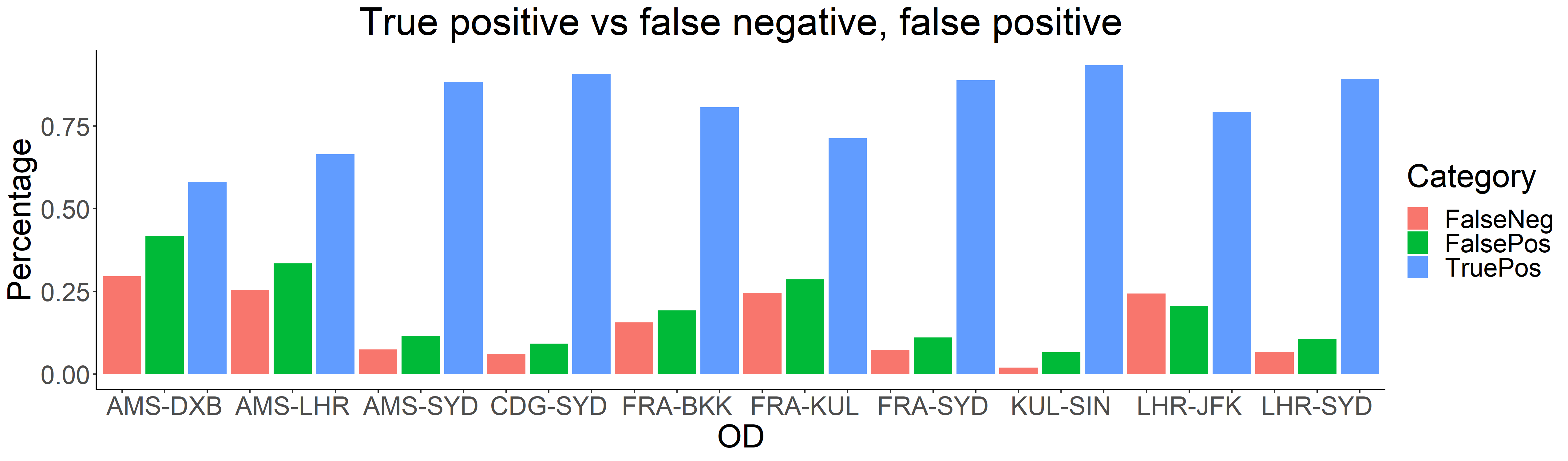} 
	\caption{Model performance by OD}
	\label{fig:posneg}
	
\end{figure}

As before, it should be noted that only 20.3\% of all options displayed were purchased. For this reason, to study the effectiveness of our model, we disregard true negatives. The percentages, in what follows, are calculated by comparing the element against true positives, false negatives and false positives. This is shown in Figure \ref{fig:posneg}.

Note the difference in performance between ODs. For example, for AMS-DXB, we note a relatively high rate of false positives and false negatives. Our initial hypothesis was this is caused by the number of competitors: the more options presented to the customer, the more dispersed the data, and therefore the more mistakes we can make. However, from Table \ref{tab:num_comp} we note that this is not the case. For example, consider the AMS-DXB with AMS-LHR ODs. These have respectively seven and four competitors. However, the number of false negatives and false positives in terms of percentage for AMS-LHR is barely different. Now compare LHR-JFK with KUL-SIN. Both of these ODs have two competitors, but the model performs much better for the latter OD. We have therefore evidence which may suggest that the LHR-JFK market is one that has more dynamics than the KUL-SIN market. 

Overall, the results across ODs are very positive. This is illustrated by the median percentage of 84.5\%. The best performing OD is KUL-SIN with a score of 93.4\%, closely followed by the SYD ODs. AMS-LHR and AMS-DXB with 66.5\% and 58.1\% respectively score worst.

Complete confusion matrices are show in the appendix. Out of all ODs, we feel that the AMS-DXB case is the most worrying. Comparing 1769 true positives with 1275 false negatives, the model severely underestimates the total number of bookings. These ratios are much lower for other ODs.

\subsection{Customer behavior}
\label{sec:booking_customerbehavior}
\begin{figure}[H]
	\includegraphics[width=1\textwidth]{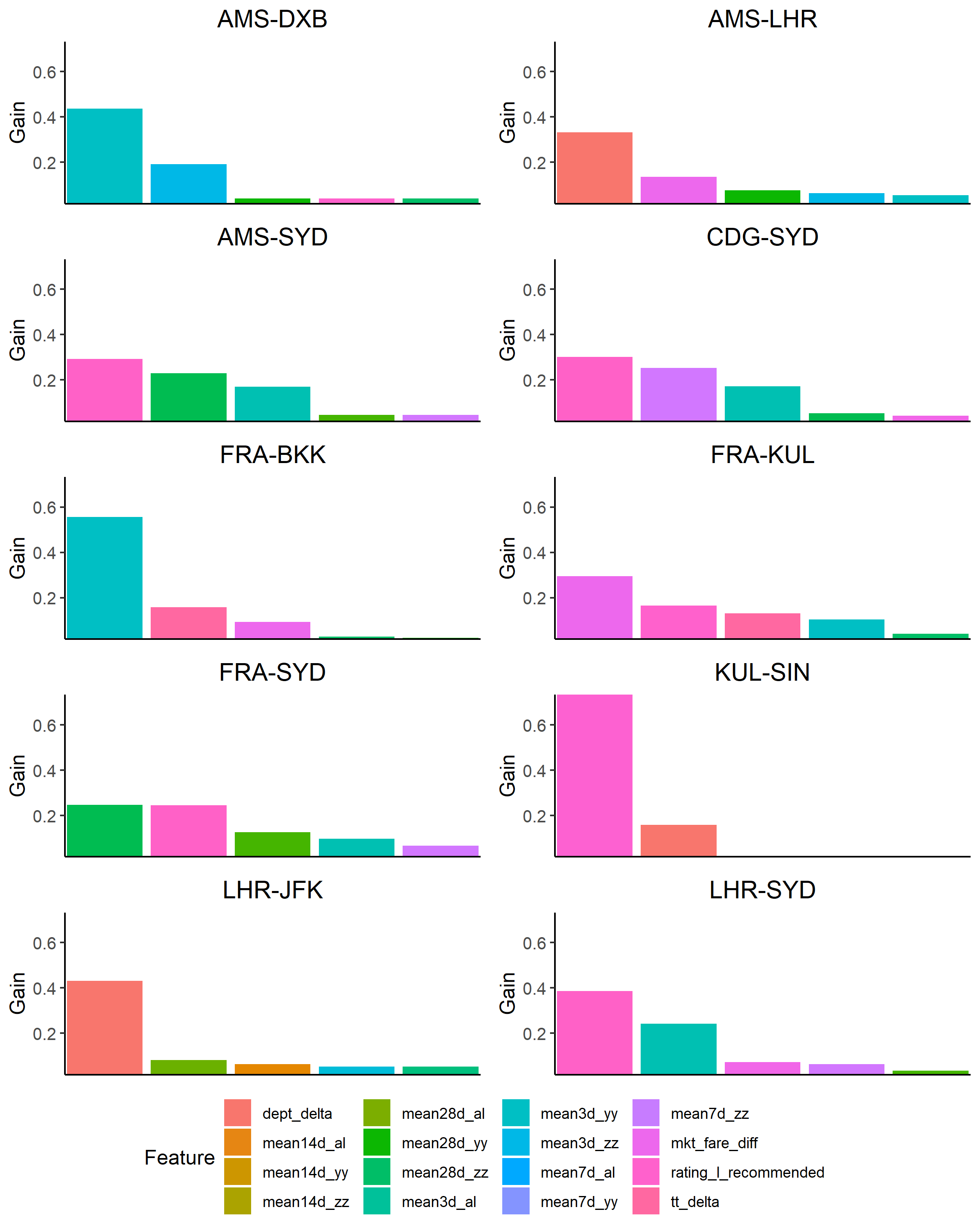} 
	\caption{Gain of features by OD, sorted by importance from left to right}
	\label{fig:importance}
\end{figure}

In our opinion, the more interesting topic is what features drive booking behavior. We refer the reader to Figure \ref{fig:importance}. This figure shows, by OD, the gain of each feature. We exclude any features with a gain smaller than $0.05$.

We note three very different behaviors: price sensitive ODs (OD pairs that are predominantly affected by price), departure-time sensitive ODs (those OD pairs that are driven by schedule) and comfort ODs (OD pairs for which passenger comfort is important). These are discussed in Sections \ref{sec:booking_results_od_pricesensitive}, \ref{sec:booking_results_od_deptime} and \ref{sec:booking_results_od_comfort} respectively.

\subsubsection{Price sensitive ODs}
\label{sec:booking_results_od_pricesensitive}
For these ODs, the pricing features dominate booking behavior. Interestingly, for AMS-DXB, the features of how the airline moves against the cheapest and second cheapest airline in the market is dominant; the actual price difference to that cheapest airline is not. Therefore, it is safe to assume that customers are set on booking with a particular airline, but monitor how this preferred airline is pricing itself against those cheaper airlines. This observation underlines the difference and added-value of our research and development of these metrics: while the price movement feature captures how the fare offering by the host airline and its competitors has moved since time of booking, the price differential feature is a snapshot of the current fare offering compared to the fare at time of booking.

Another example of this is FRA-SYD. The feature that captures the price movement against the cheapest airline in the market is dominant. On the other hand, the FRA-KUL market seems to be concerned with the fare different against the cheapest airline. The same applies for KUL-SIN – it only seems to be driven by price. To summarize, the OD's in this category are AMS-DXB, FRA-SYD, FRA-KUL and KUL-SIN.
\subsubsection{Departure time sensitive ODs}
\label{sec:booking_results_od_deptime}
Amsterdam, London and New York are traditionally considered routes for business travellers. One of the reasons for this is that all of these cities are financial hubs. It is worth noting that the delta in time between the airline’s departure time and 6AM has a lot of explanatory power. Note that especially in the LHR-JFK case, the pricing features have little weight. In the AMS-LHR case, the fare difference versus the cheapest airline in the market is the second most powerful variable. We hypothesize that this is as a result of the recent increase in low cost carrier (LCC) frequencies between Amsterdam and London. In summary, the OD's that make up this category are AMS-LHR and LHR-JFK.
\subsubsection{Flying comfort ODs}
\label{sec:booking_results_od_comfort}
On these ultra long-haul ODs, a passenger’s comfort is a deciding factor. Note that the review scores of the IFE have great explanatory power, and actually have greater explanatory power than the price. Also note that other features that may describe the comfort of a journey, such as the quality of the seat, crew or ground services did not appear.  This seems to indicate that on long-haul ODs, passengers value their entertainment more than their seat! While airlines traditionally segment their pricing based on the origin, it is worth nothing that for these ODs terminating in Sydney customer behavior seems fairly consistent. Summarizing, the OD's that fall in this category are: AMS-SYD, CDG-SYD, FRA-SYD and LHR-SYD.

\subsubsection{Comparing ODs}

It is interesting to compare the model's performance, shown in Figure \ref{fig:posneg} with the features used to obtain the predictions, shown in Figure \ref{fig:importance}.

Consider the AMS-SYD and CDG-SYD ODs. From Table \ref{tab:num_comp}, we note that these ODs have five and four competitors, respectively. Looking at Figure \ref{fig:demand_by_comp}, we note that the distribution of demand by airlines is similar: Airline 8 is missing from CDG-SYD. Airline 4 has proportionally more demand for CDG-SYD because of this missing airline. In short, from a competition perspective, we can argue that these are similar. Looking at the model's performance in Figure \ref{fig:posneg}, we note very similar results for all three metrics: false positive, false negative and true positive. The metrics that powered these predictions, in Figure \ref{fig:importance}, show that the information gain for both the rating of IFE and $mean3d_yy$ features are very similar and these combined have an information gain of $0.45$, with IFE being the most important feature. For this reason, we have segmented these together as a "comfort" OD.

\subsection{XGB performance}

\label{sec:booking_xgboost_perf}

\begin{figure}[H]
	\includegraphics[width=1\textwidth]{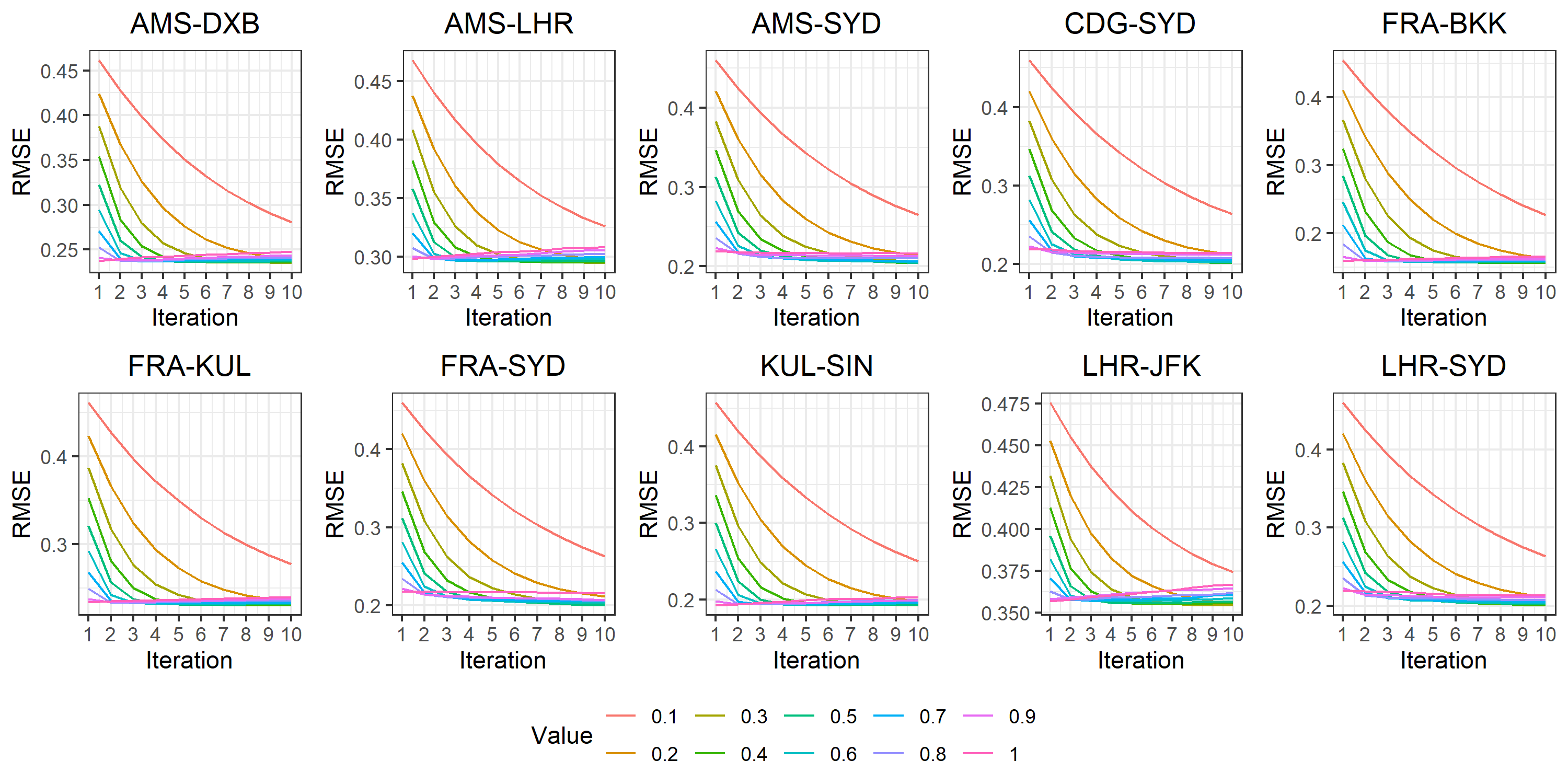} 
	\caption{Performance of learning rate $\eta$ by iteration and OD}
	\label{fig:booking_xgb_perf_eta}
\end{figure}

\begin{figure}[H]
	\includegraphics[width=1\textwidth]{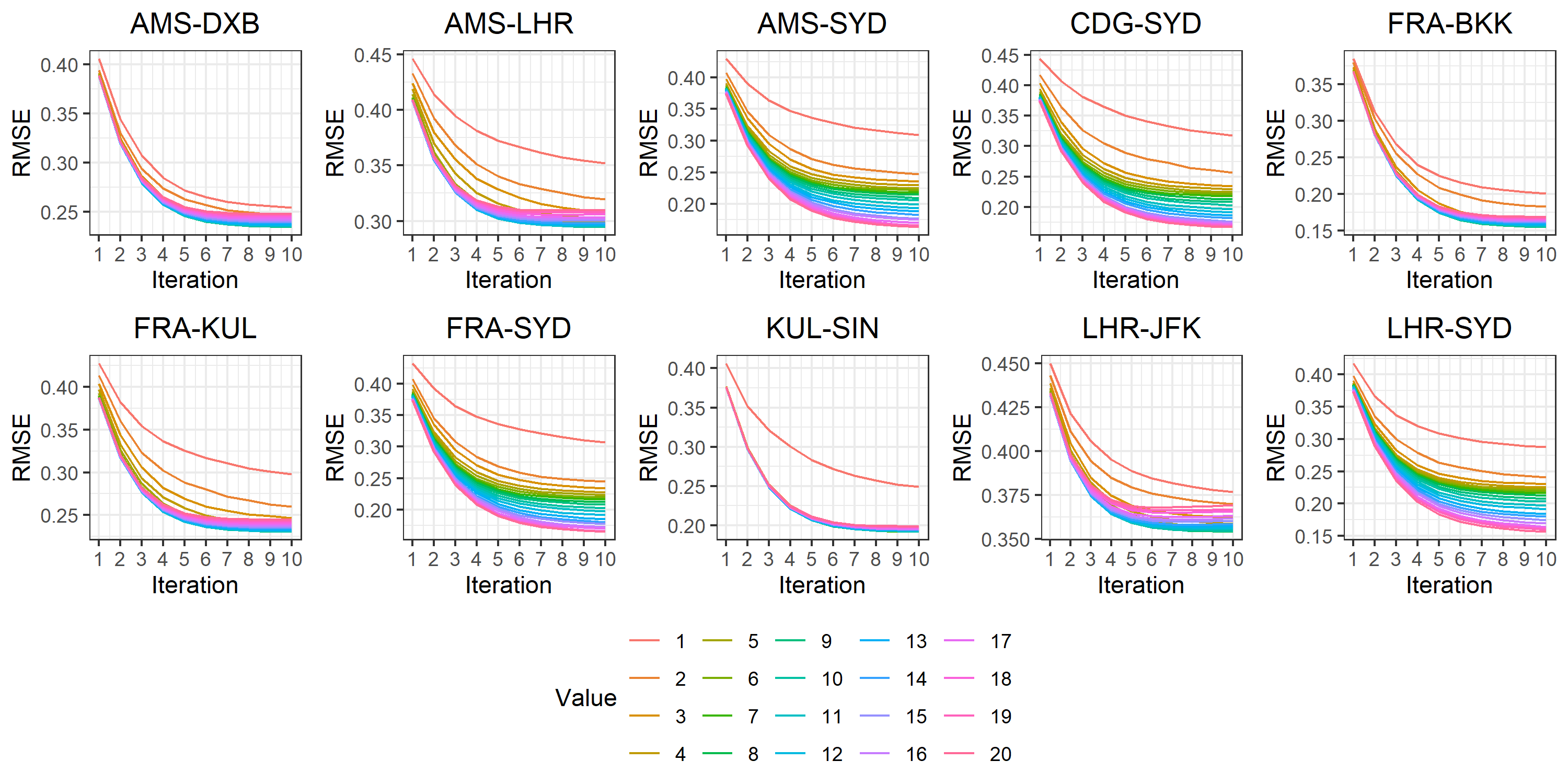} 
	\caption{Performance of $d_t$ by iteration and OD}
	\label{fig:booking_xgb_perf_treedepth}
\end{figure}

\begin{figure}[H]
	\includegraphics[width=1.1\textwidth]{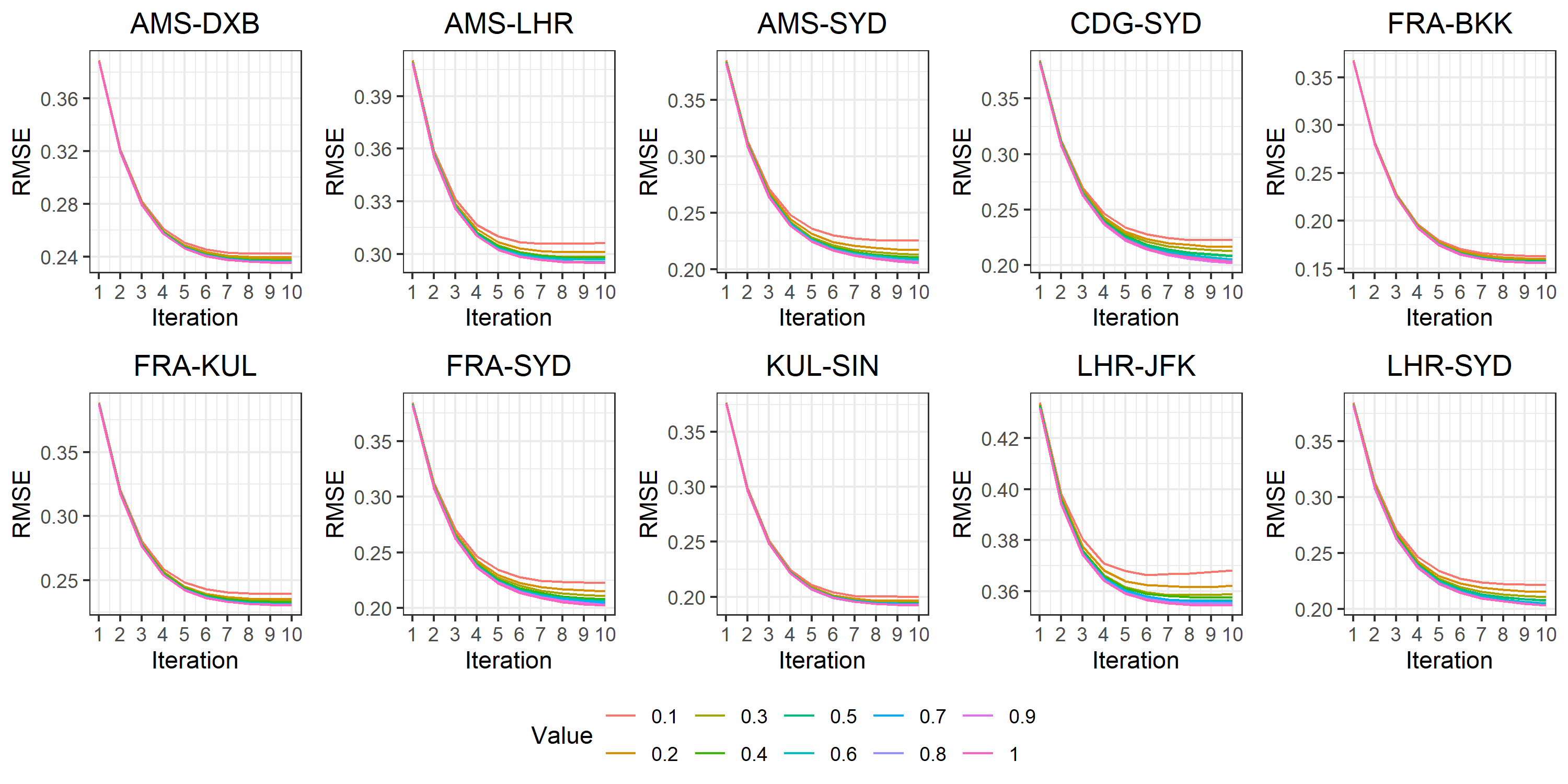} 
	\caption{Performance of $s_t$ by iteration and OD}
	\label{fig:booking_xgb_perf_sample}
\end{figure}

\begin{figure}[H]
	\includegraphics[width=1\textwidth]{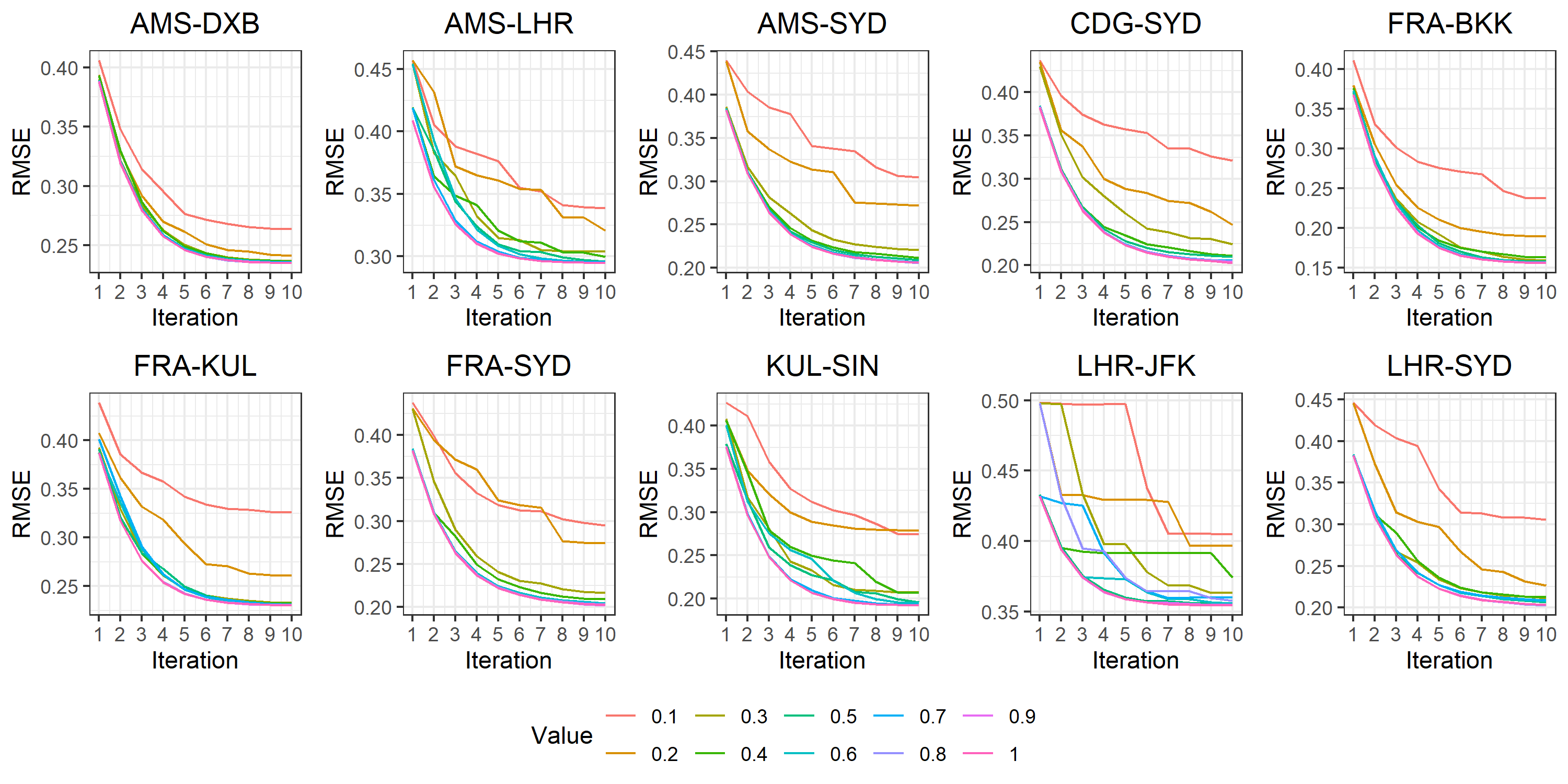} 
	\caption{Performance of $f_t$ by iteration and OD}
	\label{fig:booking_xgb_perf_colsample}
\end{figure}

In this section, we will review the parameters of the XGB model we introduced in Section \ref{sec:booking_select_params}. Note that the number of rounds is not shown here, but shown in Figure \ref{fig:error_interation} and discussed in Section \ref{sec:booking_xgb_performance}.\\
\\
\textbf{Learning Rate, $\eta$}\\
Figure \ref{fig:booking_xgb_perf_eta} shows the different values of learning rate $\eta$ for different ODs. Recall that this parameter, also known as shrinkage, controls how much weight a new tree is assigned.
From Figure \ref{fig:booking_xgb_perf_eta}, as expected, we see a very slow improvement of RMSE for $\eta = 0.1$ (red line). Increasing this parameter to $\eta = 0.2$ (orange) shows a much faster convergence. Values of around $\eta = 0.3$ (moss green) seem to be the best trade off between finding a good RMSE and runtime. For this reason, $\eta = 0.3$ is chosen. Note that very large choices of $\eta$ actually results in a slight increase of RMSE as the number of iterations grow across the different ODs. In conclusion, we do not observe great differences between the different ODs, and keep the parameter fixed at $\eta = 0.3$ for all ODs.\\
\\
\textbf{Depth of the tree, $d_t$}\\
Figure \ref{fig:booking_xgb_perf_treedepth} shows the development of RMSE on the test set for different values of the depth of the tree we allow, $d_t$.
First, let us consider the KUL-SIN OD. Here, we see a dramatic improvement from $d_t = 1$ (red), a tree with one level, to $d_t = 2$ (dark orange), a tree with two. The performance of this model, measured as the RMSE over the test set, no longer improves as we grow $d_t$ further. This seems to indicate that the ensemble of trees for this OD is relatively simple. Comparing the different ODs, we may conclude that different ODs have different levels of complexities: for the AMS-SYD, CDG-SYD, FRA-SYD and LHR-SYD ODs, we see a clear improvement as we grow $d_t$ up to its maximum (chosen) value of $d_t = 20$, while other ODs seem to converge at values of $d_t$ beyond $7$ (light blue). The value of $d_t$ depends on the OD and is chosen visually, at the lowest value after which we do not see any improvement of performance. For example, for AMS-SYD we choose $d_t = 5$ (moss green), since we do not see any improvement for larger values of $d_t$, while for LHR-SYD, we choose $d_t = 20$ (bright pink).\\
\\
\newpage
\textbf{Number of observations, $s_t$}\\
In Figure \ref{fig:booking_xgb_perf_sample}, we study the effects of the number of observations used when building trees. 
From Figure \ref{fig:booking_xgb_perf_sample}, we observe that the parameter $s_t$ is important for some ODs, while it is not important for others. For example, consider AMS-DXB. Earlier, we saw that the value of $d_t$ was important for this OD, while it does not seem to be of importance for $s_t$, even for very small values of $s_t$. This may indicate that the observations for this OD are fairly uniform: the model performs just as well randomly selecting 10\% of all observations as it does using all observations. Other ODs, such as FRA-SYD, seem to exhibit different model performance for different values. We choose the smallest value of $s_t$ accordingly: this helps reducing both runtime and avoids overfitting. For AMS-DXB, we choose $s_t = 0.2$ (orange), and for ODs that are affected by the number of observations, such as LHR-JFK, we choose $s_t = 0.1$ (red). \\
\\
\textbf{Number of features, $f_t$}\\
Figure \ref{fig:booking_xgb_perf_colsample} shows the effect of RMSE over ten iterations when using a certain proportion $f_t$ of features when constructing trees for different ODs.
In practice, a value of $f_t = 0.1$ (red) means that we use less than two features when modelling. This exemplifies the poor performance in terms of RMSE for all ODs. In Figure \ref{fig:importance} we observed that for most ODs, two features are often very important in the performance of the model. Not (randomly) choosing these two features will therefore naturally result in poor performance. While the performance of the model does improve for low values of $f_t$, the trade-off is not worth it: even when choosing $f_t = 1$, the longest runtime across the 10 different ODs is less than four seconds. In summary, we choose $f_t = 1$ for all ODs.

\subsection{Examples}

One of the main drawbacks of XGB is that method is a \textit{black box} method. This method produces an ensemble of decision trees. While a single decision tree is easy to understand, an ensemble is not. In this section, we use the \textit{xgboostExplainer} ~\cite{foster2017xgboostexplainer} package to show, for a number of examples, how the model arrives at its prediction. This is achieved by drawing the ensemble of all trees, and traversing them to obtain the probability estimation.

We will review three examples below. Each of these examples represents a given itinerary. The probability of purchase is shown on the vertical axis. The numbers inside the bar represent the log odds of each feature. As discussed in the approach section, our goal is to predict whether a itinerary is purchased. For this reason, we round the probability to the nearest integer. As a result, the horizontal line at $p=0.5$ represents the cut-off for predicting whether or not this itinerary is purchased or not. Features are ordered by their weight from left to right.

\begin{figure}[H]
	\includegraphics[width=0.85\textwidth]{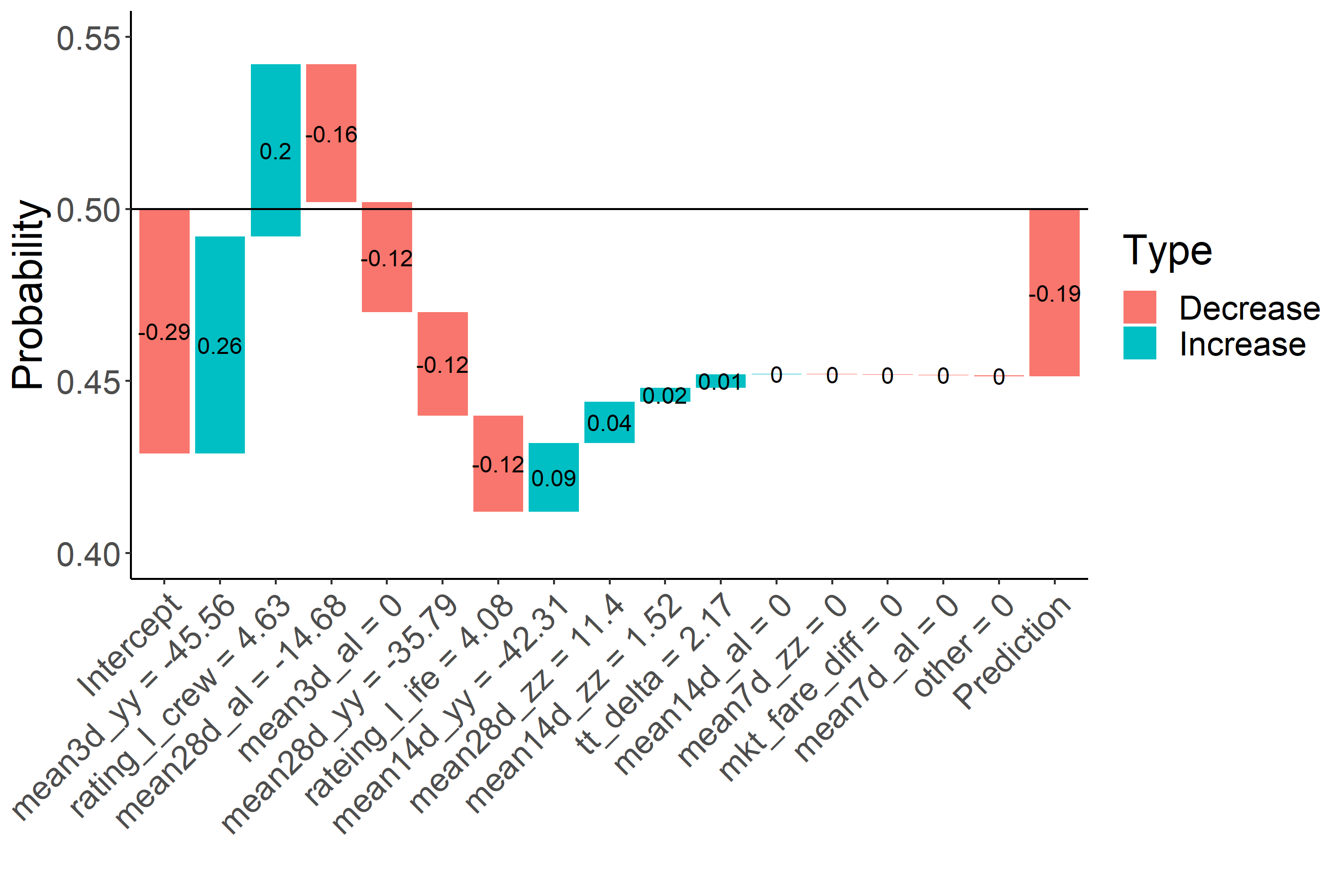} 
	\caption{Example one: long-haul itinerary, low-priced airline}
	\label{fig:behavior_ex1}
\end{figure}

Figure \ref{fig:behavior_ex1}  shows en example of how the probability is generated for a long-haul itinerary, for a low-priced airline. The model starts at a probability of purchase of 43\%, the intercept. This probability increases by 6\% to 49\%, because this airline was $\$46$ cheaper than the cheapest airline in this OD market (notice the negative sign, indicating the host airline's fare is lower is lower than the next airline in the market). The crew for this airline is rated very high, at $4.63/5$, which causes the probability of purchase to strengthen to $55\%$. This probability drops back to just over $50\%$ because the airline itself has been cheaper in the past 28 days, on average by $\$14$, indicating a price sensitive market. Note that the IFE score of $4.08$ drags the purchase probability down by almost $3\%$. The final probability of purchase is $0.45$, which, after rounding, is marked as a non-purchase.

\begin{figure}[H]
	\includegraphics[width=0.85\textwidth]{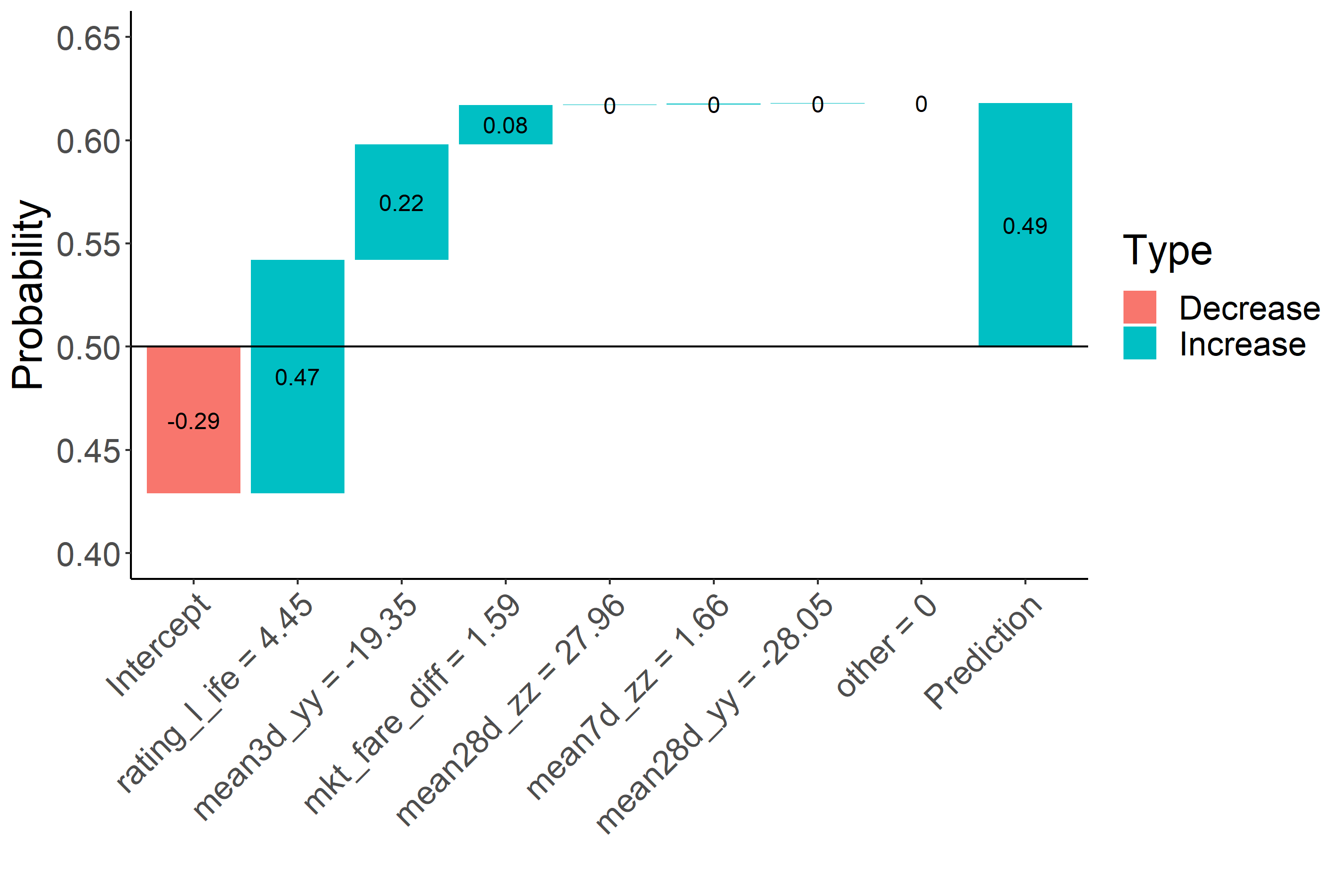} 
	\caption{Example two: long-haul itinerary, well-priced airline with great IFE}
	\label{fig:behavior_ex2}
\end{figure}

Figure \ref{fig:behavior_ex2} shows another example. In this example, the airline is well-priced compared to competition and possesses a great IFE product. The baseline purchase probability is again $43\%$. Here we clearly see the impact of IFE, the probability of purchase increases by over $11\%$. Compare this to the impact the pricing difference to the cheapest carrier in the market in the past three days:  the log odd impact of the IFE rating weighs twice as heavy as this pricing feature ($0.47$ and $0.22$, respectively). Other features have negligible impact. The model generates a purchase probability of over $62\%$.

\begin{figure}[H]
	\includegraphics[width=0.85\textwidth]{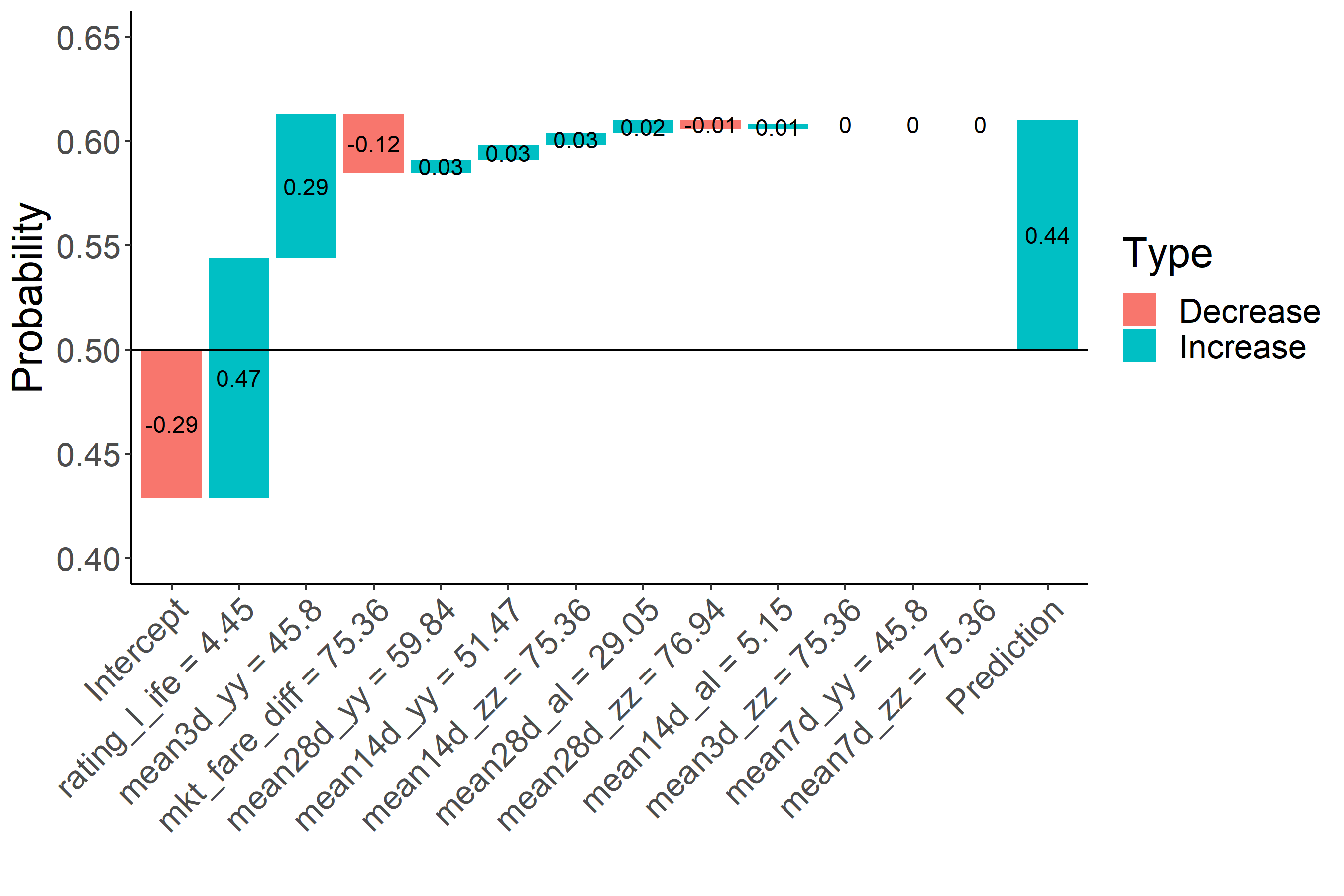} 
	\caption{Example three: long-haul itinerary, higher-priced airline with great IFE}
	\label{fig:behavior_ex3}
\end{figure}

Figure \ref{fig:behavior_ex3} shows another example of how the model generates its purchase probability. In this case, the airline has great IFE but is higher priced than competitors. The baseline purchase probability is $43\%$. Again, we see how the IFE has the biggest impact on determining purchase probability. Interestingly, the airline is on average by $\$45$ more expensive than the cheapest carrier in the market, this actually \textit{increases} the purchase probability. This could indicate the willingness to pay more for a itinerary that contains a better IFE product. However, it also shows a bound - note that because this carrier is $\$75$ more expensive than the market fare at time of purchase, the purchase probability declines slightly $2\%$. This may indicate that the sweet spot this airline could demand lies anywhere between $\$45$ and $\$75$.

\subsection{In practice}

As we discussed, predicting whether an itinerary is purchased can provide great insights in the positioning of the airline. It can benefit pricing and marketing teams and help them understand how to each of the features affect a customer's decision. While this itself provides great insights, we believe that the true value of this forecasting framework will become apparent when considering the greater revenue management problem. Consider what the end objective of RM is: maximizing revenue. Studies, for example, \cite{weatherford2002revenue}, have shown that increases in forecast accuracy show clear increases in revenue performance. Intuitively, this makes sense: inputs closer to reality should mean better output revenues.

Therefore, in this section we will compare demand forecasting through our framework with demand forecasting typically done in practice. In practice, this is done by forecasting aggregated OD/class (for brevity) demand. To compare revenues, we will use a time series forecast in simulations to provide us with a base case. Our framework, on the other hand, provides estimates on an itinerary-level. Therefore, we roll up these estimates to obtain an aggregated, OD/class-level forecast, and then use this as input in simulations for a single flight. This enables us to have a like-for-like comparison and study the effects of this method of forecasting in the end goal of revenue management: maximizing revenue.

Unfortunately, we were only given access to data for $10$ ODs. In practice, however, a flight has many different ODs crossing it. We are unable to create forecasts for every OD. Therefore, we have chosen a flight from a big airline that has four out of the ten ODs utilizing it. These four ODs combined represent, on average, $42\%$ of the total demand for this flight.

The simulation is constructed as follows:
\begin{enumerate}
    \item Calculate an optimal policy using a traditional forecasting method, in our case we chose double exponential smoothing;
    \item Calculate an optimal policy using our new forecasting framework discussed in this paper for the ODs that are available, use the traditional forecasting method for those ODs we do not have data for;
    \item Generate $n = 500$ arrival processes, with an average demand factor of $0.98$ and standard deviation of $0.1$;
    \item Simulate accepting requests and calculate revenues for both methods.
\end{enumerate}

In the simulations, we make a distinction between whether a customer downsells or not. Downsell occurs when a customer purchases a cheaper class than she was willing to pay for. Suppose we have twelve different classes for sale, with ordered fares such that $f_1$ is the most expensive product and fare $f_{12}$ the cheapest. Now consider, for example, a customer wanting to purchase booking class $10$ with fare $f_{10}$. However, class $12$ is available for sale with fare $f_{12}$. A customer downsells if she actually purchases class $12$. In this case, the airline loses $f_{10} - f_{12}$. We have done extensive research into downsell and have covered this in the literature \cite{downsell_paper}.

The fares for the four ODs we have used are covered in Table \ref{tab:booking_fares}. The associated demand by fare brand is given in Table \ref{tab:booking_demand}. Demand by booking class has a relatively high variance, so we have summarized these in terms of fare brand. Recall that the definition of a fare brand was introduced in \cite{ourpaper_lit}: fare brands are a collection of products with identical fare conditions, with the only differing factor price. 

\begin{table}[H]
\centering
\resizebox{\textwidth}{!}{
\begin{tabular}{l|cccccccccccc}
  \hline
 OD \textbackslash Class & 1 & 2 & 3 & 4 & 5 & 6 & 7 & 8 & 9 & 10 & 11 & 12 \\ 
  \hline
 1 & 2324 & 1913 & 1672 & 1152 & 1081 & 966 & 871 & 706 & 660 & 498 & 494 & 447 \\ 
  2 & 2489 & 2078 & 1995 & 1707 & 1462 & 1363 & 1187 & 1009 & 774 & 553 & 534 & 474 \\ 
  3 & 1904 & 1621 & 1323 & 1094 & 1091 & 962 & 922 & 737 & 682 & 622 & 495 & 311 \\ 
  4 & 2509 & 2043 & 1452 & 1420 & 1035 & 762 & 700 & 523 & 449 & 374 & 311 & 206 \\ 
   \hline
\end{tabular}}
\caption{Fare by booking class}
\label{tab:booking_fares}
\end{table}

\begin{table}[H]
\centering
\begin{tabular}{l|ccc}
  \hline
 OD \textbackslash Fare Brand & 1 & 2 & 3 \\ 
  \hline
1 & 0.05 & 0.30 & 0.65 \\ 
  2 & 0.10 & 0.18 & 0.70 \\ 
  3 & 0.02 & 0.46 & 0.40 \\ 
  4 & 0.12 & 0.20 & 0.59 \\ 
   \hline
\end{tabular}
\caption{Percentage of demand by fare brand. Booking Classes $1$ through $3$ are part of Fare Brand 1, $4$ through $8$ part of Fare Brand 2, $9$ through $12$ part of Fare Brand 3.}
\label{tab:booking_demand}
\end{table}

Table \ref{tab:booking_rev_perf} shows the simulation results.

\begin{table}[h]
\centering
\begin{tabular}{l|rrr}
  \hline
 Downsell & Std & XGB & \% Gain \\ 
  \hline
 No & 27468 & 27648 & 0.70 \\ 
  Yes & 9622 & 9898 & 2.90 \\ 
   \hline
\end{tabular}
\caption{Revenue performance comparing forecast methods Std, Standard forecasting, with XGB. Demand factor = $0.98$}
\label{tab:booking_rev_perf}
\end{table}

In Table \ref{tab:booking_rev_perf}, we compare the standard forecasting technique (Std) with XGB. The numbers represent the revenue performance for this flight. The percentage gain shows the relative performance. If customers do not downsell, our simulations show an average improvement of 0.7\%. However, of particular interest is the scenario where customers do downsell. This is a more realistic scenario, and this provides revenue gains of 2.90 \%. We therefore have found evidence that this framework is particularly beneficial in these cases.

\section{Discussion}
\label{sec:booking_discussion}

The results show our ability to segment ODs in three distinct categories (price-, schedule- and comfort-sensitive). Currently, we use bookings from the entire period for which we have been given access to, as one big dataset. It could be that there is seasonality in what customers prefer. For example, in the summer peak, when schools have their holidays, people may be less price sensitive. After all, they are bound by these dates to travel, much more so than in an off-peak period, when it is easy to move your travels a day or two.

Another detail we would like to highlight is that the data shown here is from one specific OTA. This brings us to the following. Firstly, it could be that there are airlines that do not sell ticket through this OTA. For example, in the US, Southwest Airlines only sells tickets through their own channel, not through any OTA. Secondly, while we dealt with a large enough number of bookings to make our research statistically significant, one has to remember that the airline also sells through other channels, including other OTAs. For this reason, the dataset may not represent the entire population. However, we cannot think of reasons why customer behavior on one OTA will be drastically different than behavior on another channel.

When using ratings from airlines, we chose to calculate aggregates by airline, not by OD and airline. As discussed, we took this decision because the number of reviews would get extremely small and as a result aggregates are not reliable. Particularly for airlines that offer many different types of aircraft - some of them having certain features, such as WiFi or IFE, while other aircraft do not - this may not be fair to aggregate review scores like this.

While elements such as safety records may be a deciding factor intuitively, the number of accidents or incidents is at an historically low level. In effect, flying has become very safe, and therefore, there is little difference in measures between airlines. For this reason, we suspect the model is not using this feature for any OD.

Intra-day changes to fares can cause the fares from the booking dataset mismatch those from the competitive pricing dataset. While bookings depend on (real time) availability, the fares in the competitive pricing dataset are only scraped once a day (but always at the same time).

The implications for RM are as follows. We have shown that on some ODs, passengers only seem to look at fare in their decision-making process. As a result, the airline could considering lowering fares for these ODs. Alternatively, it could look into unbundling fares to become more competitive. On the other hand, the approach in this paper extends to other departments of airlines: it shows that the investments some airlines make in their entertainment offering do drive extra bookings. 

In traditional RM optimisation techniques, one optimizes network revenue using demand forecasts, estimated fares as well as capacity constraints. Some techniques assume booking class-independent demand forecasts, while others use customer choice probabilities to derive demand forecasts. In this paper, we have presented a framework to estimate whether a given itinerary will get purchased or not. Naturally, rolling up these estimates from an itinerary to OD-level, we can derive a true, competitor-based demand forecast.

Another implication for RM is having an ability to estimate what your airline's fare premium should be, given your offering. It is often a discussion in airlines how much more expensive, or similarly, how much more cheaper you should be compared to your competitors given your product. For figure \ref{fig:behavior_ex3} we estimated what this airline's premium could be for this OD, given its very well-rated IFE offering.

Before discussing the results, we compared the logit to the XGB model. We found that the XGB model outperforms the logit model for true positives. However, we also found that the logit model is more biased toward predicting false negatives while the XGB model is more biased toward false positives. In practice, we would prefer a false positive to some extent. After all, predicting a false negative is arguably worse than a false positive: research typically shows that overforecasting results in better revenue performance than underforecasting. It may be interesting to study why these methods are biased toward different errors, but this is outside the scope of this paper.

We argue that the most important result of our work is obtaining a better forecast. Weatherford et al \cite{weatherford2002revenue}, show that overforecasting by 25\% accuracy lead to a revenue loss of 1.35\%, while overforecasting by 12.5\% lead to revenues 0.18\% lower than an optimal forecast. The fares ranged from \$66 to \$275 in their work, with demand consistent across classes. Comparing our work to this work, we suspect that our results are higher as a result of a wider-spread fare ladder and a different optimization technique. Since our forecast considers competition, this could be a straightforward method to incorporate competitor-based information into the optimization process. One should remember that the results presented in this paper, in particular in Table \ref{tab:booking_rev_perf}, are a lower bound to actual revenue performance. After all, only 42\% of the flight's total demand was modelled using our new technique. The objective of this paper was to introduce a new method to predict booking behavior. The results we showed have a demand factor of $0.98$. We suggest further research into different scenarios - for example, extremely empty of extremely popular flights; different optimization techniques; different fare ladders to study the effects of our method on revenue performance.

Finally, the method used in this paper is a black box method. As we discussed, while a decision tree is easy to understand, an ensemble of decision trees is not. Therefore, the analyst may have its doubts on how the model arrives at its prediction. The \textit{xgboostExplainer} tool, which we used to produce Figures \ref{fig:behavior_ex1}, \ref{fig:behavior_ex2} and \ref{fig:behavior_ex3} is a great tool to give insight in the model, which in turn will restore analyst confidence.

\section{Appendix: engineered features}
\label{sec:booking_appendix}

\begin{table}[H]
\centering
\resizebox{1\textwidth}{!}{
\begin{tabular}{l|lll}
Field Title & Field Description & Source & Calculation \\
	\hline
	airline.x & Airline & Competitor Pricing &  \\
	od    & OD    & Competitor Pricing &  \\
	airline\_num\_id & Obfuscated Airline ID & Schedule &  \\
	num   & Obfuscated date & OTA   &  \\
	t.x   & Days before departure & OTA   &  \\
	home\_carrier & Is this airline a home carrier? & Schedule & Is OD's origin or destination the airline's hub? \\
	sc1   & Review score site 1 & Review Website &  \\
	sc2   & Review score site 2 & Review Website &  \\
	price & Price of itinerary & Competitor Pricing &  \\
	rating\_l\_recommended & Airline recommended by Leisure Passengers & Review Website & Median rating of all review \\
	rating\_l\_review & Review sentiment by Leisure Passengers & Review Website & Text mining based on AFINN dataset \\
	rating\_l\_fb & F\&B rating by Leisure Passengers & Review Website & Median rating of all review \\
	rating\_l\_ground & Ground services rating by Leisure Passengers & Review Website & Median rating of all review \\
	rating\_l\_ife & IFE rating by Leisure Passengers & Review Website & Median rating of all review \\
	rating\_l\_crew & Crew rating by Leisure Passengers & Review Website & Median rating of all review \\
	rating\_l\_seat & Seat rating by Leisure Passengers & Review Website & Median rating of all review \\
	rating\_l\_value & Value for money rating by Leisure Passengers & Review Website & Median rating of all review \\
	rating\_l\_wifi & WiFi rating by Leisure Passengers & Review Website & Median rating of all review \\
	rating\_l\_obs & Number of review observations by Leisure Passengers & Review Website & Count of number of reviews \\
	rating\_b\_recommended & Airline recommended by Business Passengers & Review Website & Median rating of all review \\
	rating\_b\_review & Review sentiment by Business Passengers & Review Website & Text mining based on AFINN dataset \\
	rating\_b\_fb & F\&B rating by Business Passengers & Review Website & Median rating of all review \\
	rating\_b\_ground & Ground services rating by Business Passengers & Review Website & Median rating of all review \\
	rating\_b\_ife & IFE rating by Business Passengers & Review Website & Median rating of all review \\
	rating\_b\_crew & Crew rating by Business Passengers & Review Website & Median rating of all review \\
	rating\_b\_seat & Seat rating by Business Passengers & Review Website & Median rating of all review \\
	rating\_b\_value & Value for money rating by Business Passengers & Review Website & Median rating of all review \\
	rating\_b\_wifi & WiFi rating by Business Passengers & Review Website & Median rating of all review \\
	rating\_b\_obs & Number of review observations by Business Passengers & Review Website & Count of number of reviews \\
	sent\_mean & Mean sentiment score & Review Website & Median rating of all review \\
	sent\_sd & Standard deviation sentiment score & Review Website & Text mining based on AFINN dataset \\
	sent\_mean\_rel\_diff & Difference to mean sentiment score & Review Website & Median rating of all review \\
	sent\_mean\_rel\_perc & Percentage difference to mean sentiment score & Review Website & Median rating of all review \\
	sent\_sd\_rel\_diff & Difference to sd sentiment score & Review Website & Median rating of all review \\
	sent\_sd\_rel\_perc & Percentage difference to sd sentiment score & Review Website & Median rating of all review \\
	direct\_flight & Is this a direct flight yes/no? & Review Website & Median rating of all review \\
	has\_night\_flight & Does this airline offer a night flight? & Review Website & Median rating of all review \\
	has\_day\_flight & Does this airline offer a day flight? & Review Website & Median rating of all review \\
	first\_flight\_dep & Airline's time of departure of first flight of the day & Review Website & Count of number of reviews \\
	first\_flight\_arr & Airline's time of arrival of first flight of the day & Review Website & Median rating of all review \\
	last\_flight\_dep & Airline's time of departure of last flight of the day & Review Website & Text mining based on AFINN dataset \\
	last\_flight\_arr & Airline's time of arrival of lasst flight of the day & Review Website & Median rating of all review \\
	min\_flying\_time & Airline's minimum flying time & Review Website & Median rating of all review \\
	min\_conn\_time & Airine's minimum connection time & Review Website & Median rating of all review \\
	min\_travel\_time & Airline's minimum travel time & Review Website & Median rating of all review \\
	has\_night\_departure & Is this a night departure? & Review Website & Median rating of all review \\
	has\_morning\_arrival & Is this a morning arrival? & Review Website & Median rating of all review \\
	num\_frequencies & Number of frequencies offered by airline & Review Website & Median rating of all review \\
	aircraft\_type & Aircraft type & Review Website & Count of number of reviews \\
	airline\_fleet\_size & Airline fleet size & Kaggle & Sum of airframes \\
	airline\_fleet\_cost & Airline fleet cost (estimated) & Kaggle & Sum of airframe cost \\
	airline\_fleet\_age & Airline fleet age & Kaggle & Average age of airframes \\
	bucket\_t & Bucketed time (time before departure grouped in multiples of 10) & OTA   & floor(Time before departure/10)*10 \\
  
\end{tabular}}
\end{table}

\begin{table}[H]
\centering
\resizebox{1\textwidth}{!}{
\begin{tabular}{l|lll}
Field Title & Field Description & Source & Calculation \\
	\hline
	mean3d\_al & Mean last 3 day price moment of own airline & Competitor Pricing & Mean of fare at t.x - fare at t.y, rolling 3 days, airline vs itself \\
	min3d\_al & Min of last 3 day price moment of own airline & Competitor Pricing & Min of fare at t.x - fare at t.y, rolling 3 days, airline vs itself \\
	max3d\_al & Max of last 3 day price moment of own airline & Competitor Pricing & Max of fare at t.x - fare at t.y, rolling 3 days, airline vs itself \\
	sd3d\_al & SD of last 3 day price moment of own airline & Competitor Pricing & SD of fare at t.x - fare at t.y, rolling 3 days, airline vs itself \\
	mean7d\_al & SD of last 7 day price moment of own airline & Competitor Pricing & Mean of fare at t.x - fare at t.y, rolling 7 days, airline vs itself \\
	min7d\_al & SD of last 7 day price moment of own airline & Competitor Pricing & Min of fare at t.x - fare at t.y, rolling 7 days, airline vs itself \\
	max7d\_al & SD of last 7 day price moment of own airline & Competitor Pricing & Max of fare at t.x - fare at t.y, rolling 7 days, airline vs itself \\
	sd7d\_al & SD of last 7 day price moment of own airline & Competitor Pricing & SD of fare at t.x - fare at t.y, rolling 7 days, airline vs itself \\
	mean14d\_al & SD of last 14 day price moment of own airline & Competitor Pricing & Mean of fare at t.x - fare at t.y, rolling 14 days, airline vs itself \\
	min14d\_al & SD of last 14 day price moment of own airline & Competitor Pricing & Min of fare at t.x - fare at t.y, rolling 14 days, airline vs itself \\
	max14d\_al & SD of last 14 day price moment of own airline & Competitor Pricing & Max of fare at t.x - fare at t.y, rolling 14 days, airline vs itself \\
	sd14d\_al & SD of last 14 day price moment of own airline & Competitor Pricing & SD of fare at t.x - fare at t.y, rolling 14 days, airline vs itself \\
	mean28d\_al & SD of last 28 day price moment of own airline & Competitor Pricing & Mean of fare at t.x - fare at t.y, rolling 28 days, airline vs itself \\
	min28d\_al & SD of last 28 day price moment of own airline & Competitor Pricing & Min of fare at t.x - fare at t.y, rolling 28 days, airline vs itself \\
	max28d\_al & SD of last 28 day price moment of own airline & Competitor Pricing & Max of fare at t.x - fare at t.y, rolling 28 days, airline vs itself \\
	sd28d\_al & SD of last 28 day price moment of own airline & Competitor Pricing & SD of fare at t.x - fare at t.y, rolling 28 days, airline vs itself \\
	mean3d\_yy & Mean last 3 day price moment of cheapest airline & Competitor Pricing & Mean of fare at t.x - fare at t.y, rolling 3 days, airline vs cheapest airline \\
	min3d\_yy & Min of last 3 day price moment of cheapest airline & Competitor Pricing & Min of fare at t.x - fare at t.y, rolling 3 days, airline vs cheapest airline \\
	max3d\_yy & Max of last 3 day price moment of cheapest airline & Competitor Pricing & Max of fare at t.x - fare at t.y, rolling 3 days, airline vs cheapest airline \\
	sd3d\_yy & SD of last 3 day price moment of cheapest airline & Competitor Pricing & SD of fare at t.x - fare at t.y, rolling 3 days, airline vs cheapest airline \\
	mean7d\_yy & SD of last 7 day price moment of cheapest airline & Competitor Pricing & Mean of fare at t.x - fare at t.y, rolling 7 days, airline vs cheapest airline \\
	min7d\_yy & SD of last 7 day price moment of cheapest airline & Competitor Pricing & Min of fare at t.x - fare at t.y, rolling 7 days, airline vs cheapest airline \\
	max7d\_yy & SD of last 7 day price moment of cheapest airline & Competitor Pricing & Max of fare at t.x - fare at t.y, rolling 7 days, airline vs cheapest airline \\
	sd7d\_yy & SD of last 7 day price moment of cheapest airline & Competitor Pricing & SD of fare at t.x - fare at t.y, rolling 7 days, airline vs cheapest airline \\
	mean14d\_yy & SD of last 14 day price moment of cheapest airline & Competitor Pricing & Mean of fare at t.x - fare at t.y, rolling 14 days, airline vs cheapest airline \\
	min14d\_yy & SD of last 14 day price moment of cheapest airline & Competitor Pricing & Min of fare at t.x - fare at t.y, rolling 14 days, airline vs cheapest airline \\
	max14d\_yy & SD of last 14 day price moment of cheapest airline & Competitor Pricing & Max of fare at t.x - fare at t.y, rolling 14 days, airline vs cheapest airline \\
	sd14d\_yy & SD of last 14 day price moment of cheapest airline & Competitor Pricing & SD of fare at t.x - fare at t.y, rolling 14 days, airline vs cheapest airline \\
	mean28d\_yy & SD of last 28 day price moment of cheapest airline & Competitor Pricing & Mean of fare at t.x - fare at t.y, rolling 28 days, airline vs cheapest airline \\
	min28d\_yy & SD of last 28 day price moment of cheapest airline & Competitor Pricing & Min of fare at t.x - fare at t.y, rolling 28 days, airline vs cheapest airline \\
	max28d\_yy & SD of last 28 day price moment of cheapest airline & Competitor Pricing & Max of fare at t.x - fare at t.y, rolling 28 days, airline vs cheapest airline \\
	sd28d\_yy & SD of last 28 day price moment of cheapest airline & Competitor Pricing & SD of fare at t.x - fare at t.y, rolling 28 days, airline vs cheapest airline \\
	mean3d\_zz & Mean last 3 day price moment of second cheapest airline & Competitor Pricing & Mean of fare at t.x - fare at t.y, rolling 3 days, airline vs second cheapest airline \\
	min3d\_zz & Min of last 3 day price moment of second cheapest airline & Competitor Pricing & Min of fare at t.x - fare at t.y, rolling 3 days, airline vs second cheapest airline \\
	max3d\_zz & Max of last 3 day price moment of second cheapest airline & Competitor Pricing & Max of fare at t.x - fare at t.y, rolling 3 days, airline vs second cheapest airline \\
	sd3d\_zz & SD of last 3 day price moment of second cheapest airline & Competitor Pricing & SD of fare at t.x - fare at t.y, rolling 3 days, airline vs second cheapest airline \\
	mean7d\_zz & SD of last 7 day price moment of second cheapest airline & Competitor Pricing & Mean of fare at t.x - fare at t.y, rolling 7 days, airline vs second cheapest airline \\
	min7d\_zz & SD of last 7 day price moment of second cheapest airline & Competitor Pricing & Min of fare at t.x - fare at t.y, rolling 7 days, airline vs second cheapest airline \\
	max7d\_zz & SD of last 7 day price moment of second cheapest airline & Competitor Pricing & Max of fare at t.x - fare at t.y, rolling 7 days, airline vs second cheapest airline \\
	sd7d\_zz & SD of last 7 day price moment of second cheapest airline & Competitor Pricing & SD of fare at t.x - fare at t.y, rolling 7 days, airline vs second cheapest airline \\
	mean14d\_zz & SD of last 14 day price moment of second cheapest airline & Competitor Pricing & Mean of fare at t.x - fare at t.y, rolling 14 days, airline vs second cheapest airline \\
	min14d\_zz & SD of last 14 day price moment of second cheapest airline & Competitor Pricing & Min of fare at t.x - fare at t.y, rolling 14 days, airline vs second cheapest airline \\
	max14d\_zz & SD of last 14 day price moment of second cheapest airline & Competitor Pricing & Max of fare at t.x - fare at t.y, rolling 14 days, airline vs second cheapest airline \\
	sd14d\_zz & SD of last 14 day price moment of second cheapest airline & Competitor Pricing & SD of fare at t.x - fare at t.y, rolling 14 days, airline vs second cheapest airline \\
	mean28d\_zz & SD of last 28 day price moment of second cheapest airline & Competitor Pricing & Mean of fare at t.x - fare at t.y, rolling 28 days, airline vs second cheapest airline \\
	min28d\_zz & SD of last 28 day price moment of second cheapest airline & Competitor Pricing & Min of fare at t.x - fare at t.y, rolling 28 days, airline vs second cheapest airline \\
	max28d\_zz & SD of last 28 day price moment of second cheapest airline & Competitor Pricing & Max of fare at t.x - fare at t.y, rolling 28 days, airline vs second cheapest airline \\
	sd28d\_zz & SD of last 28 day price moment of second cheapest airline & Competitor Pricing & SD of fare at t.x - fare at t.y, rolling 28 days, airline vs second cheapest airline \\
	dep\_time\_mam & Departure time in minutes after midnight & Schedule & Hours and minutes converted into minutes \\
	connecting\_time & Connecting time & Schedule & Dep time next flight - Arr time previous flight \\
	travel\_time & Travel time & Schedule & Flying time + connecting time \\
	mintt & Minimum connecting time offered by airline & Schedule & Min(travel\_time) by OD by days before departure \\
	tt\_delta & Difference between airline travel time and min travel time for this OD & Schedule & Airline's travel time - mintt \\
	dept\_delta & Difference between itinerary's departure time and "ideal" deparure time & Schedule & Departure time - 7AM departure \\
	mkt\_fare & Lowest fare in the market & Competitor Pricing & Min fare by OD by time before departure \\
	mkt\_fare\_diff & Difference to lowest fare in the market & Competitor Pricing & Airline's fare - mkt\_fare \\
	mkt\_fare\_diff\_perc & Percentage difference to lowest fare in the market & Competitor Pricing & (Airline's fare/mkt\_fare)-1 \\
	is\_cheapest & Is this airline cheapest? & Competitor Pricing & If airline's fare = mkt\_fare \\
	is\_bought & Itineray purchased? Label & OTA   & Is this itinerary bought? \\
	airline\_id & Obfuscated airline ID & OTA   & Randon number \\

\end{tabular}}
\caption{Overview of all engineered features}
\label{tab:features}
\end{table}

\bibliographystyle{plain}
\bibliography{bibs} 

\end{document}